\documentclass[twoside,a4paper,10pt]{amsart}

\usepackage{graphicx}

\usepackage{amsfonts,amssymb,amsmath}
\usepackage[all,cmtip]{xy}

\newcommand{\MV}{\mathrm{MV}_\infty}
\newcommand{\eL}{\mathbb{L}}
\newcommand{\Cm}{\mathcal{C}}
\newcommand{\Pp}{\mathcal{P}}
\newcommand{\Qm}{\mathcal{Q}}
\newcommand{\caL}{\mathcal{L}}

\renewcommand{\epsilon}{\varepsilon}

\newcommand{\ac}{\overline{\mathrm{ac}}}
\newcommand{\ord}{\mathrm{ord}}

\newcommand{\Rr}{\mathbb{R}}
\newcommand{\Nn}{\mathbb{N}}
\newcommand{\Zz}{\mathbb{Z}}
\newcommand{\Cc}{\mathbb{C}}
\newcommand{\Qq}{\mathbb{Q}}

\newcommand{\Aa}{\mathbb{A}}

\newcommand{\Dd}{\mathcal{D}}
\newcommand{\Ppr}{\mathbb{P}}
\newcommand{\mrm}[1]{\;\mathrm{#1}\;}

\newcommand{\Jac}{\mathrm{Jac}}

\newcommand{\set}[1]{\left\{#1\right\}}

\newcommand{\abs}[1]{\left| #1 \right|}

\newcommand{\card}[1]{\mathrm{card}\left( #1 \right)}
\newcommand{\lcm}{\mathrm{lcm}}

\newcommand{\ie}{\emph{i.e.}}

\newcommand{\kk}{\mathbf{k}}

\newcommand{\fonction}[5]{#1 :\left\{\begin{array}{ccl}
#2&\to&#3\\[3mm]
#4&\mapsto&#5
\end{array}\right.}

\theoremstyle{plain}
\newtheorem{nth}{theorem}[section]

\newtheorem{theorem}[nth]{Theorem}
\newtheorem{proposition}[nth]{Proposition}
\newtheorem{lemma}[nth]{Lemma}
\newtheorem{corollary}[nth]{Corollary}

\theoremstyle{definition}
\newtheorem{definition}[nth]{Definition}
\newtheorem{lemma-Definition}[nth]{Lemma-Definition}
\newtheorem{proposition-definition}[nth]{Proposition-Definition}

\newtheorem{example}[nth]{Example}

\theoremstyle{remark}
\newtheorem{remark}[nth]{Remark}

\begin{document}

\title{Motivic local density}
\author{Arthur Forey}
\address{{A. Forey, Institut de Math\'ematiques de Jussieu, UMR 7586 du CNRS, Universit\'e Pierre et Marie Curie, 4 place Jussieu, 75252 Paris Cedex 05}}
\email{arthur.forey@imj-prg.fr}
\thanks{This research was partially supported by ERC grant 246903 NMNAG and ANR-15-CE40-0008 (D\'efig\'eo).}
\date{\today}


\begin{abstract}
We develop a theory of local densities and tangent cones in a motivic framework, extending work by Cluckers-Comte-Loeser about $p$-adic local density. We prove some results about geometry of definable sets in Henselian valued fields of characteristic zero, both in semi-algebraic and subanalytic languages, and study Lipschitz continuous maps between such sets. We prove existence of regular stratifications satisfying analogous of Verdier condition $(w_f)$. Using Cluckers-Loeser theory of motivic integration, we define a notion of motivic local density with values in the Grothendieck ring of the theory of the residue sorts. We then prove the existence of a distinguished tangent cone and that one can compute the local density on this cone endowed with appropriate motivic multiplicities. As an application, we prove a uniformity theorem for $p$-adic local density.

\end{abstract}

\maketitle

\section{Introduction}

Notions of local densities have been developed in various contexts since their introduction by Lelong in \cite{lelong_integration_1957}. We recall what is known in the complex and real cases. Denote by $\lambda_d$ the Lebesgue measure on $\Rr^n$ of dimension $d$ and $B^n(a,r)$ the ball around $a\in \Rr^n$ of radius $r$. If $X\subseteq \Rr^n$ of dimension $d$ and $a\in X$, the local density of $X$ at $a$ is the limit, if it makes sense, 
\[
\Theta_d(X,a)=\lim_{r\to 0} \frac{\lambda_d(X\cap B^n(a,r))}{\lambda_d(B^d(0,r))}.
\]
In \cite{lelong_integration_1957}, Lelong first shows the existence of the limit if $X$ is a complex analytic manifold of real dimension $d$. In this complex case, Thie in \cite{thie_lelong_1967} expresses the local density as a sum of the local densities of components of the tangent cone of $X$ at $a$ with multiplicities, which implies that $\Theta_d(X,a)$ is a positive integer. Then in \cite{draper_intersection_1969}, Draper interprets the local density as the algebraic multiplicity of the local ring of $X$ at $a$. 

In \cite{kurdyka_densite_1989}, Kurdyka and Raby study the real subanalytic case and show the existence of the local density if $X$ is real subanalytic. They prove an analogous of Thie's result by expressing it as a finite sum of densities of components of the positive tangent cone, with multiplicities, but the local density is no longer an integer. Lion gives in \cite{lion_densite_1998} a geometric interpretation of the local density using the Cauchy-Crofton formula and proves existence of the density of semi-pfaffian sets. Comte, Lion and Rolin show in \cite{comte_nature_2000} that the local density viewed as a function of the base point is a log-analytic function. In fact Kurdyka and Raby approach works in any real o-minimal setting, see for example \cite{comte_equisingularite_2008} for a proof of the existence of local density in this setting.

Cluckers, Comte and Loeser study in \cite{cluckers_local_2012} the local density in a $p$-adic context. They show the existence of a notion of local density for a semi-algebraic or subanalytic subset $X\subseteq K^n$, where $K$ is a finite extension of the field of $p$-adic numbers $\Qq_p$. The definition is more subtle than just the analogous of the real one, with Lebesgue measure replaced by Haar measure, because the limit does not exist in general. However, one can consider the following sequence of normalized volumes
\[
\theta_n=\frac{\mu_{d,K}(X\cap B(a,n))}{\mu_{d,K}(B(0,n))},
\]
where $\mu_{d,K}$ is the Haar measure on $K^d$, extended to $K^n$ and 
\[B(a,n)=\set{x\in K^n\mid \ord(x-a)\geq n}\]
 is the ball around $a$ of valuative radius $n$. Then they show that there exists an integer $e>0$ such that for any $i=0,...,e-1$, the sequence $(\theta_{i+ke})_{k\geq 0}$ converges to some $d_i$. Then they define the local density of $X$ at $a$ as the mean value
\[
\Theta_d^K(X,a)=\frac{1}{e}\sum_{i=0}^{e-1} d_i.
\]
In this $p$-adic setting, there is no natural notion of tangent cone, therefore they are led to study for any $\Lambda$, a multiplicative subgroup of $K^\times$ of finite index, a $\Lambda$-tangent cone. They show that for a given $X$, for $\Lambda$ small enough, all the $\Lambda$-tangent cones coincide, hence there exists a distinguished tangent cone. Then they assign multiplicities to the tangent cone and show as in the real case that the local density can be computed on this tangent cone with multiplicities. One should point out that in \cite{cluckers_local_2012}, two different notions of multiplicities are introduced and claimed to coincide but it is not the case, see Section \ref{cor-padic} ; only the most naive one using cardinality of fibers over the tangent cone leads to the correct notion of multiplicities. 

In this paper, we generalize this notion. Fix $L=\Cc((t))$ the field of formal Laurent series over the complex numbers. In this introduction, we limit ourself to this field for simplicity, but in the next sections we work in greater generality, considering some Henselian fields of characteristic zero and all residue characteristic. We also do not work with a single fixed field but with a tame or mixed tame theory of fields, see Section \ref{section_tame_theory}. Consider the three sorted language of Denef-Pas, with one sort for the valued field $L$, with the language $\set{0,1,+,-,\cdot,t}$, one sort for the residue field $k_L=\Cc$ with the language $\set{0,1,+,-,\cdot}$ and one sort for the value group $\Zz$ with the Presburger language $\set{0,+,\leq, \set{\equiv_n}_{n\geq1}}$, where $\equiv_n$ is a binary relation for congruence modulo $n$. We also add a function $\ord : L^\times \to \Zz$ which is interpreted as the valuation on $K^\times$ and a function $\ac : L^\times \to k_L$ which we interpret by the angular component
\[
\ac\left(\sum_{i\geq i_0} a_i t^i\right)=a_{i_0} \mrm{if} a_{i_0}\neq 0.
\]
In this introduction, by definable we mean definable in this first order language. Consider $K_0(\mathrm{Var}_\Cc)$, the Grothendieck group of varieties over $\Cc$. By variety, we mean a separated scheme of finite type over $\Cc$ but a more naive definition would lead to the same group. Recall that $K_0(\mathrm{Var}_\Cc)$ is the abelian group generated by symbols $[X]$ for $X$ a variety with relations 
\[
[X]=[Y]+[X\backslash Y],
\]
where $Y$ is a closed subvariety of $X$. The cartesian product induces a product on $K_0(\mathrm{Var}_\Cc)$, making it a commutative ring (with unity the class of a point). Denote $\eL=[\Aa_\Cc^1]$ the class of the affine line and consider the localization 
\[
\mathcal{M}_\Cc=K_0(\mathrm{Var}_\Cc)\left[\eL^{-1}, \set{\frac{1}{1-\eL^{-i}}}_{i>0}\right].
\]
In \cite{cluckers_constructible_2008}, Cluckers and Loeser assign to some definable subsets $X$ of $L^n$ an element of $\mathcal{M}_\Cc$, called the motivic volume of $X$. They define more generally a whole class of so-called motivic constructible functions that are stable under integration and with good properties like Fubini's theorem, change of variable formula and integration with respect to parameters.  Denote $\mu_d$ the $d$-dimensional measure defined in \cite{cluckers_constructible_2008}, and let $B^n(a,m)=\set{x\in L^n\mid \ord(x-a)\geq m}$, the $n$-dimensional ball around $a\in L^n$ of valuative radius $m\in \Zz$. Fix $X\subseteq L^n$ definable of dimension $d$, $a\in L^n$ and define the sequence of normalized volumes by 
\[
\theta_m=\frac{\mu_d(X\cap B^n(a,m))}{\mu_d(B^d(0,m))}\in\mathcal{M}_\Cc.
\]
We prove that there exists an integer $e>0$ such that for any $i=0,1,...,e-1$, the subsequence $(\theta_{i+ke})_{k\geq 0}$ converges to some $d_i\in\mathcal{M}_\Cc$ for the dimension topology. We then define the motivic local density of $X$ at $a$ by 
\[
\Theta_d(X,a)=\frac{1}{e}\sum_{i=0}^{e-1}d_i\in \Qq\otimes\mathcal{M}_\Cc.
\]
This periodic convergence can occur even if $X$ is an algebraic set, for example the set $X=\set{x^2=y^3}$ has a density of 1/2. 

In order compute the local density of the tangent cone, we have to consider $\Lambda$-tangent cones $C_a^\Lambda(X)$ for all $\Lambda$ for the form 
\[\Lambda_n=\set{\lambda\in L^\times\mid \ac(\lambda)=1,\ord(\lambda)\equiv 0\mod n}.\]
We are able to show that there is a distinguished tangent cone, that is, there is a $\Lambda$ such that for any smaller $\Lambda'$, $C_a^{\Lambda'}(X)=C_a^{\Lambda}(X)$. This distinguished tangent cone captures the local properties of $X$ at $a$, in the sense that one can compute the local density of $X$ at $a$ on this cone. To do so, one needs to assign multiplicities to this cone. Whereas in the $p$-adic case the multiplicities are integers, here the multiplicities need to be some motives. In order to define the multiplicities, we partition $X$ into graphs of 1-Lipschitz functions parametrized by (a definable of) the residue field. This allows us to define the tangent cone with multiplicities which is a motivic constructible function of support $C_a^{\Lambda}(X)$. 
Our main Theorem \ref{thm-princ} is an analogue of respective theorems of Thie, Kurdyka-Raby and Cluckers-Comte-Loeser, stating that the local density of $X$ is the local density of its tangent cone with multiplicities. By the specialization principle of Cluckers-Loeser, this implies a uniformity result for $p$-adic local density.

In Section \ref{section_preliminaries}, we define the framework of (mixed-)tame theories, prove some results about Lipschitz continuous function, a weak curve selection lemma and recall the setting of Cluckers-Loeser theory of motivic integration. In Section \ref{section_local_density}, we define the motivic local density and the tangent cone with multiplicities in order to state our main theorem. Section \ref{section_regular_stratifications} is devoted to establishing the existence of regular stratification satisfying analogies of Verdier condition $(w_f)$. Section \ref{section_proof_main_theorem} is devoted to the proof of the main theorem. In Section \ref{section_application_padic}, we first show how our results specialize to the $p$-adic case and then prove a uniformity theorem for $p$-adic local density.

\vspace{1cm}
I would like to thank my advisor Fran\c{c}ois Loeser for suggesting me to work on this subject and many useful pieces of advice. I would also like to thank Raf Cluckers and Georges Comte for useful discussions. I also thank the referee for the numerous comments.

\section{Preliminaries}
\label{section_preliminaries}
\subsection{Tame theories of valued fields}

\label{section_tame_theory}

A discretely valued field $K$ is a field with a surjective group homomorphism
\[
\ord : K^\times \to \Zz,
\]
extended by $\ord(0)=\infty$ and satisfying for all $x,y\in K^\times$, 
\[
\ord(x+y)\geq \min\set{\ord(x),\ord(y)}.
\]
A ball in $K$ is a set of the form $\set{x\in K \mid \ord(x-a)\geq \alpha}$ for some $a\in K$ and $\alpha \in \Zz$. The collection of all balls of $K$ forms a basis of the so-called valuation topology on $K$. The field $K$ is called Henselian if the valuation ring $\mathcal{O}_K$ is Henselian. Write $\mathfrak{M}_K$ for the maximal ideal of $\mathcal{O}_K$. 

\begin{definition}
Fix $p$ to be either 0 or a prime number, and an integer $e$. A $(0,p,e)$-field is a Henselian discretely valued field $K$ of characteristic 0, of residue characteristic $p$, together with a chosen uniformizer $\varpi_K$ of the valuation ring $\mathcal{O}_K$ of $K$, such that either $0=p=e$ or $p>0$ and the ramification degree is $e$, that is, $\ord(\varpi_K^e)=\ord(p)=e$. 
\end{definition}

We will always identify the value group of a $(0,p,e)$-field with $\Zz$. For example, $\Qq_p$ with $p$ as uniformizer is a $(0,p,1)$-field. A $(0,p,e)$-field comes with higher order angular components maps, defined for $n\geq 1$, 
\[
\ac_n : K^\times \to R_{n,K}:=\mathcal{O}_K/\mathfrak{M}_K^n,\; x \mapsto \varpi_K^{-\ord(x)} x \mod \mathfrak{M}_K^n,
\]
extended by $\ac_n(0)=0$. Write $\ac=\ac_1$, note that the maps $\ac_n$ are multiplicative on $K$ and coincide on $\mathcal{O}_K$ with the natural projection $\mathcal{O}_K \to R_{n,K}$. Note also that $R_{1,K}=k_K$ is the residue field of $K$. We call the $R_{n,K}$ the residue rings of $K$. 

\begin{definition}[First order language]
We define now the semi-algebraic language that we will use to describe subsets of $(0,p,e)$-fields. Let $\mathcal{L}_{\mathrm{high}}$ the first order language consisting of one sort for the valued field, one sort for the value group, and one sort for each residue ring. On these sorts, $\mathcal{L}_{\mathrm{high}}$ consists of the language of rings $(0,1,+,-,\cdot)$ and a symbol $\pi$ for the uniformizer on the valued field, the language of Presburger arithmetics $(0,1,+,-,\leq,\set{\cdot\equiv\cdot \mod n}_{n> 1})$ on the value group, the language of rings on each residue ring, a symbol $\ord$ for the valuation map, symbols $\ac_n$ for the higher order angular components, and for $n\geq m$, the symbol $p_{n,m}$ for natural projection between residue rings $R_n$ and $R_m$. 
We will also use the Denef-Pas language $\mathcal{L}_{\mathrm{DP}}$ which is the restriction of $\mathcal{L}_{\mathrm{high}}$ with only three sorts : one for the valued field, one for the value group and one for the residue field and only one angular component $\ac$. 
\end{definition}

Any $(0,p,e)$-field has a natural interpretation as $\mathcal{L}_{\mathrm{high}}$-structure. Let $\mathcal{T}_{(0,p,e)}$ the $\mathcal{L}_{\mathrm{high}}$-theory of sentences true in every $(0,p,e)$-field. As concrete list of axioms that implies this theory can be found in the work of Pas \cite{pas_cell_1990}, in particular, $\mathcal{T}_{(0,p,e)}$ eliminates field quantifiers. 

Let $\mathcal{L}$ a language with the same sorts as $\mathcal{L}_{\mathrm{high}}$, containing $\mathcal{L}_{\mathrm{high}}$, and $\mathcal{T}$ be an $\mathcal{L}$-theory containing $\mathcal{T}_{(0,p,e)}$. 

The following axioms will be about properties of every model $\mathcal{K}$ of $\mathcal{T}$, even those which are not $(0,p,e)$-fields (this happens if the value group is not isomorphic to $\Zz$). By abuse of notation, we will often identify $\mathcal{K}$, a model of $\mathcal{T}$, with the associated valued field $K$. 

We call the sort for the valued field the main sort, and the other ones the auxiliary sorts. We call the collection of all the sorts for the residue rings the residue sorts. 

We will need the following notion of $\kk$-partition, in order to work piecewise over the residue rings. 
\begin{definition}[Definable $\kk$-partition]
If $X$ is a definable set of a model $K$ of $\mathcal{T}$, a definable $\kk$-partition of $X$ is an injective definable function
\[
g : X \to \prod_{i=1}^s R_{n_i,K},
\]
for some integers $n_i\geq 1$. A fiber of $g$ is called a $\kk$-part of $X$. 

Given two definable $\kk$-partitions $g$ and $g'$ of $X$, $g'$ is said to be a refinement of $g$ if every $\kk$-part of $X$ for $g'$ is included into a $\kk$-part of $X$ for $g$.  It is clear that given two definable $\kk$-partitions, we can find a third one which is a refinement of the two others. 
\end{definition}

If $X\subseteq A\times B$, and $b\in B$, we denote the fiber of $X$ over $b$ by
\[ X_b := \set{a\in A\mid (a,b)\in X}.\]

\begin{definition}[Jacobian property for functions]
Let $K$ be a model of $\mathcal{T}_{(0,p,e)}$. Let $F : X\subseteq K \to K$ be a function. We say that $F$ has the Jacobian property on $X$ if for every ball $B\subseteq X$, either $F$ is constant on $B$, or the following properties hold :
\begin{enumerate}
\item $F$ is injective on $B$ and $F(B)$ is a ball,
\item $F$ is $C^1$ on $B$ with derivative $F'$,
\item $F'$ is nonzero on $B$ and $\ord(F')$ is constant on $B$,
\item for all $x,y\in B$ with $x\neq y$, we have
\[
\ord(F')+\ord(x-y)=\ord(F(x)-F(y)).
\]
\end{enumerate}
Moreover, if $n>0$ is an integer, we say that $F$ has the $n$-Jacobian property if the following additional conditions hold :
\begin{enumerate}
\setcounter{enumi}{5}
\item $\ac_n(F')$ is constant on $B$,
\item for all $x,y\in B$ with $x\neq y$, we have
\[
\ac_n(F')\cdot \ac_n(x-y)=\ac_n(F(x)-F(y)).
\]
\end{enumerate}
\end{definition}

\begin{definition}[Jacobian property for $\mathcal{T}$]
Say that the Jacobian property holds for the $\mathcal{L}$-theory $\mathcal{T}$ if any model $K$ satisfies the following. For any integer $n\geq 1$, any definable function $F : X\subset K\times Y \to K$, there is a definable $\kk$-partition of $X$ such that for any $y\in Y$ and any $\kk$-part of $X_y$, the restriction of $F$ to this set has the Jacobian property. In more explicit terms, if $g$ is the function witnessing the $\kk$-partition, then for every $\xi$ and every $y\in Y$, then
\[
x\in g^{-1}(\xi)_y \mapsto F(x,y)
\]
has the Jacobian property.
\end{definition}

The following definition translates the stable embedding of the value group and orthogonality between value group and residue rings. 
\begin{definition}[Split]
Say that the theory $\mathcal{T}$ splits if the following holds. For any model $\mathcal{K}$, with valued field $K$, value group $\Gamma$, and residue rings $R_{n,K}$, we have :
\begin{enumerate}
\item any $\mathcal{K}$-definable subset of $\Gamma^r$ is $\Gamma$-definable in the language $(+,-,0,<)$,
\item any definable subset $X\subseteq \Gamma^r\times\prod_{i=1}^s R_{n_i,K}$ is equal to a finite disjoint union of $Y_i\times Z_i$ where the $Y_i$ are definable subsets of $\Gamma^r$ and the $Z_i$ are definable subsets of $\prod_{i=1}^s R_{n_i,K}$. 
\end{enumerate}
\end{definition}

The notion of $b$-minimality has been introduced in \cite{cluckers_b-minimality_2007}, it is closely related to the concept of cell decomposition. It is shown in \cite{cluckers_motivic_????}, Corollary 3.8 that the following notion of finite $b$-minimality implies $b$-minimality.
\begin{definition}[Finite $b$-minimality]
Call $\mathcal{T}$ finitely $b$-minimal if the following holds for any model $K$. Each locally constant definable function $g : K \to K$ has finite image, and for any definable set $X\subseteq K\times Y \to K$, where $Y\subseteq K^n$ is definable, there exists an integer $\ell\geq 1$, a $\kk$-partition of $X$
\[
g : X \to A,
\]
a definable function
\[
c : Y\times A \to K,
\]
and a definable $B\subseteq \Gamma\times Y$
such that for every $\xi \in A$ and $y\in Y$ such that if $g^{-1}(\xi)_y$ is nonempty, $g^{-1}(\xi)_y$ is either
\begin{enumerate}
\item equal to the singleton $\set{c(y,\xi)}$,
\item equal to the ball $\set{x\in K \mid \ac_\ell(x-c(y,\xi))=\zeta, \ord(x-c(y,\xi))\in B_y}$ where $\zeta$ is one of the component of $\xi$. 
\end{enumerate}
In the first case, we say $X$ is a 0-cell, in the second we say it is a 1-cell. We call $c$ the center of the cell, although it is not necessarily unique and does not need to lie in the cell. We say that the set $g^{-1}(\xi)_y$ is in $c(y,\xi)$-config. 
\end{definition}

The abstract notion of tame theories is introduced by Cluckers and Loeser in \cite{cluckers_motivic_????}, but examples of such theories were already previously considered, see Example \ref{example_tame} below. 
 
\begin{definition}[Tame theory]
Let $\mathcal{T}$ a theory containing $\mathcal{T}_{(0,p,e)}$ in a language $\mathcal{L}$ with the same sorts as $\mathcal{L}_{\mathrm{high}}$, which splits, is finitely $b$-minimal, has the Jacobian property and has at least one $(0,p,e)$-field as a model. Then we call $\mathcal{T}$ a tame theory of fields if $e=0$ and a mixed tame theory of fields if $e>0$. 

If $\mathcal{T}$ is an $\mathcal{L}$-theory with he same sorts as $\mathcal{L}_{\mathrm{DP}}$, contains $\mathcal{T}_{(0,0,0)}$, splits, is finitely $b$-minimal and has the Jacobian property, we also call $\mathcal{T}$ a tame theory of fields.
A $\mathcal{T}$-field is a model of $\mathcal{T}$ that is a $(0,p,e)$-field. 
\end{definition}

\begin{example}
\label{example_tame}
\begin{enumerate}
\item The first examples of (mixed) tame theory are given by Pas in \cite{pas_uniform_1989} and \cite{pas_cell_1990} : the theory $\mathcal{T}_{(0,p,e)}$ is a (mixed) tame theory of fields. This is the semi-algebraic example. 

\item \label{example_analytic}
We now define an analytic language. Let $K$ a $(0,p,e)$-field. Let $K\{x_1,...,x_n\}$ the ring of formal power series $\sum_{i\in \Nn} a_i x^i$ over $K$ such that $\ord(a_i)$ goes to $+\infty$ when the $i$ goes to $+\infty$. If $K$ is complete, this coincide with power series converging on $\mathcal{O}_K^n$. Then let $\mathcal{L}_{K,\mathrm{an}}$ the language $\mathcal{L}_{\mathrm{high}}$ with a function symbol for each element of $K\{x_1,...,x_n\}$ (for all $n>0$). Then each $(0,p,e)$-complete extension $L$ of $K$ has a natural 
$\mathcal{L}_{K,\mathrm{an}}$-structure, a function $f\in K\{x_1,...,x_n\}$ is interpreted as the corresponding function $\mathcal{O}_L^n\to L$ and extended by 0 outside $\mathcal{O}_L^n$. Let $\mathcal{T}_{K,\mathrm{an}}$ the theory of all such fields $L$. Then $\mathcal{T}_{K,\mathrm{an}}$  is a (mixed) tame theory by \cite{CLR_analytic_2006}. 

\item More generally, for every analytic structure defined in \cite{cluckers_fields_2011}, with associated language $\mathcal{L}$ and theory $\mathcal{T}$, if $\mathcal{T}$ has a $(0,p,e)$ model, then $\mathcal{T}$ is a (mixed) tame theory.

\item Let $\mathcal{T}$ a (mixed) tame theory in a language $\mathcal{L}$, and $R$ be a subring of a $\mathcal{T}$-field. Let $\mathcal{L}(R)$ the language $\mathcal{L}$ with additional constant symbols for elements of $R$. Any $\mathcal{T}$-field that is an algebra over $R$ and with $\ord$ and $\ac_n$ extending the one on $R$, can be viewed as an $\mathcal{L}(R)$-structure. Let $\mathcal{T}(R)$ the $\mathcal{L}(R)$-theory of such structures. Then $\mathcal{T}(R)$ is a (mixed) tame theory, by \cite{cluckers_motivic_????}, Proposition 3.11. 
\end{enumerate}
\end{example}

\subsection{Definable subassignments}
Fix a (mixed-)tame theory of fields $\mathcal{T}$. Fix $n,\ell,s$ some nonnegative integers, and $r=r_1,\dots,r_\ell$ a tuple of $\ell$ nonnegative integers. Define $h[n,r,s]$ to be the functor form the category of $\mathcal{T}$-fields to the category of sets sending any $\mathcal{T}$-field $K$ to 
\[
h[n,r,s](K)=K^n\times R_{r_1,K}\times \cdots \times R_{r_\ell,K}\times \Zz^s.
\]
A definable subassignment of $h[n,r,s]$ is the data for every $\mathcal{T}$-field $K$ of a subset $X_K$ of $h[n,r,s](K)$ such that for some $\mathcal{L}$-formula $\varphi$ (with free variables of corresponding lengths), $X_K=\varphi(K)$. Most often, we will just say that $X$ is definable and call $X_K$ a definable subset. 
A definable function $f : X \to Y$ is a definable subassignment $G$ such that $G(K)$ is the graph of a function $f_K : X_K\to Y_K$. 

As usual in model theory, we fix a $\mathcal{T}$-field $K$ with residue rings sufficiently saturated, such that the study of definable subassignments $X$ in $\mathcal{T}$ is reduced to the study of definable subsets $X_K$

\subsection{Dimension theory}

A (mixed-)tame theory is $b$-minimal in the sense of \cite{cluckers_b-minimality_2007}, see \cite[Corollary 3.8]{cluckers_motivic_????}. In particular, one can use dimension theory for definable in $b$-minimal theories. Fix a (mixed-)tame theory $\mathcal{T}$ and a $\mathcal{T}$-field $K$. 
The dimension of a definable set $X\subseteq K^n$ is the biggest $d$ such that for some coordinate projection $p : K^n \to K^d$, $p(X)$ contains an open ball of $K^d$ (with convention $\dim(\emptyset)=-\infty$). 

If $X\subseteq K^n\times S$, where $S$ is a product of residue rings and $\Zz^s$, then the dimension of $X$ is the dimension of $p(X)$ with $p$ the projection to $K^n$. 

If $f : X \to Y$ is definable, then the dimension of $X$ over $Y$ (also called dimension of $X$ relatively to $Y$) is the maximum of the dimensions of the fibers $f^{-1}(y)$ over $y\in Y$. Denote $\dim(X)$ the dimension of $X$ and $\dim(X/Y)$ the dimension of $X$ over $Y$ (which depends implicitly on $f$).  

If $X$ is a definable subassigment, it may happen that $\dim(X_K)\neq \dim(X_{L})$ for $K$ and $L$ two $\mathcal{T}$-fields. For example, for $X=\set{x\mid x^2+1=0}$,  $K=\Rr((t))$, $L=\Cc((t))$, then $\dim(X_K)=-\infty$ and $\dim(X_L)=0$.

Hence we define $\dim(X)$ to be the maximum of the $\dim(X_K)$ for $K$ a $\mathcal{T}$-field. However, if $K$ has its residue rings sufficiently saturated, then $\dim(X)=\dim(X_K)$, hence we will ignore this subtlety in the rest of the paper. 

Here are some properties of the dimension that we will use thorough the paper. 
\begin{proposition}[{\cite[Propositions 4.2 and 4.3]{cluckers_b-minimality_2007}}]

Let $W,Y,Z$ be definable, $X,X',C$ be definable over $Y$, with $X,X'\subseteq C$, and $f : X\to X'$ definable over $Y$ (that is, commuting to projections to $Y$). Then
\begin{itemize}
\item $\dim(A\cup B/Y)=\max\set{\dim(A/Y),\dim(B/Y)}$,
\item $\dim(W\times Z)=\dim(W)+\dim(Z)$,
\item $\dim(X)\geq \dim(f(X)$,
\item For each integer $d\geq 0$, the set $S_f(d)=\set{x'\in X'\mid \dim(f^{-1}(x')/Y)=d)}$ is definable and
\[
\dim(f^{-1}(S_f(d))/Y)=\dim(S_f(d)/Y)+d,
\]
with convention $-\infty+d=-\infty$,
\item if $Y$ is contained in some auxiliary sorts, then $\dim(X/Y)=\dim(X)$. 
\end{itemize}
\end{proposition}

\subsection{}

We fix a (mixed-)tame theory of fields $\mathcal{T}$ and $K$ a $\mathcal{T}$-field. 

We will use either additive or multiplicative notation for the value group, denoting respectively $\ord(x)$ and $\abs{x}=e^{-\ord(x)}$ the valuation and the absolute value of an element $x\in K$. It is a purely conventional matter; depending on the context one might be more readable than the other. 
If $x=(x_i)_{i}\in K^n$, we set $\ord(x)=\min_i{\ord(x_i)}$, this defines a topology on $K^n$. If $X\subseteq K^n$, $\overline{X}$ the closure of $X$, $\mathring{X}$ the interior of $X$ and $\partial X$ the boundary of $X$ are all definable.
Define $B(a,m)=\set{x\in K^n\mid \ord(x-a)\geq m}$ to be the ball around $a$ of valuative radius $m$. 

\begin{definition}
Consider a function $f : A \to B$. Then $f$ is said globally $C$-Lipschitz if for all $x,y\in A$, 
\[
\abs{f(x)-f(y)}\leq C \abs{x-y}.
\]

The function $f$ is locally $C$-Lipschitz at $u\in A$ if there is an $n\in\Zz$ such that $f$ is globally $C$-Lipschitz on $U\cap B(u,n)$. 

The function $f$ is locally $C$-Lipschitz if it is locally $C$-Lipschitz at all $u\in A$. 
\end{definition}

In particular, if $f$ is differentiable, then $f$ is locally 1-Lipschitz if $\ord(\partial f/\partial x)\geq 0$.

\begin{proposition}[Injectivity versus constancy {\cite[2.2.5]{cluckers_non-archimedean_????}}]
\label{inj-cons}
Let $F : X\subseteq K \to K$ be definable. Then there exists a $\kk$-partition of $X$ such that the restriction of $F$ on each $\kk$-part is either injective or constant. 
\end{proposition}

The following property is used in \cite{cluckers_motivic_????} to define motivic integral. In positive residue characteristic, one can already see the need of using higher angular components. 
\begin{lemma}
\label{lemma-finite-bij}
Let $X\subseteq K^n$ be finite and definable. Then $X$ is in definable bijection with a subset of some residue rings.
\end{lemma}
\begin{proof}
 We work by induction on the number $n$ of points. If $X=\set{x_1,...,x_n}$, set $y=\frac{1}{n} \sum_{i=1}^n x_i$, then there are $e>0$ and $i\neq j$ such that $\ac_e(x_i-y)\neq \ac_e(x_j-y)$, then this reduces to the case of less than $n$ points.
\end{proof}

The following proposition will be useful in the study of the tangent cone, allowing us consider only tangent cones of graphs of 1-Lipschitz functions, which are easier to study.

\begin{proposition}
\label{partition-1lipsch}
Let $X\subseteq K^{n}\times S$  be definable of dimension $d$ relatively to $S$. Then there exists a definable $Y\subseteq X$ of dimension $<d$ relatively to $S$, and a definable $\kk$-partition of $X\setminus Y$ such that for any $s\in S$, on each $\kk$-part $X_{s,\xi}$, there is an $(s,\xi)$-definable open $U_{s,\xi}\subseteq K^d$ and a $(s,\xi)$-definable function $\phi_{s,\xi} : U_{s,\xi} \to K^{n-d}$ of graph $\Gamma_{s,\xi}$ that is differentiable and locally 1-Lipschitz and such that for a permutation of coordinates $\gamma_{s,\xi}$, $X_{s,\xi}=\gamma_{s,\xi}(\Gamma_{s,\xi})$.
\end{proposition}

\begin{proof} 
Working piecewise in the fibers of the projection of $X\to S$ and using logical compactness, one can suppose $X\subseteq K^n$. One proceeds by induction on $n$. If $n=d$ one proceeds $Y=X\setminus \mathring{X}$, and we are done. If $n>d$, for $i=1,..., n$, one considers the projection $\pi_i : K^{n} \to K^{n-1}$ that omits the $i$-th coordinate. Then $X$ is partitioned into two definable parts $X_0^i$ and $X_1^i$, such that $\pi_{i \vert X}$ has 0-dimensional fibers on $X_0^i$ and 1-dimensional fibers on $X_1^i$. The set $\pi_i(X_1^i)\subseteq K^{n-1}$ is then of dimension $d-1$, so one is done by the induction hypothesis. On $X_0^i$, the fibers are 0-dimensional, so they are in definable bijection with a definable in the auxiliary sorts. Up to adding some $\kk$-parameters, one can assume the fibers are one point, $i.e.$ $\pi_i : X_0^i \to \pi_i(X_0^i)$ is a $\kk$-definable bijection. One is then reduced to considering the case $X=\cap_i X_0^i$, \emph{i.e.} $\pi_i : X \to \pi_i(X)$ is a $\kk$-definable bijection for all $i=1,...,n$. Denote $h_i$ the $i$-th coordinate function of $\pi_i^{-1}$. By the Jacobian property, up to adding some additional $\kk$-parameters and neglecting a set of dimension $<d$, one can assume all $h_i$ have partial derivatives on the whole of $\pi_i(X)$. Up to cutting $X$ into finitely many pieces, one can also assume that the order of the partial derivatives of the functions $h_i$ is either nonnegative or negative on the whole $\pi_i(X)$. It suffices to find an $i_0$ such that all the partial derivatives of $h_{i_0}$ are nonnegative and then we are done by induction. But using the chain rule for derivation, this is clear that such an $i_0$ exists. 
 \end{proof}

\subsection{Some results about Lipschitz continuity}

We recall now some results of \cite{cluckers_non-archimedean_????} about reaching global Lipschitz continuity  from local Lipschitz continuity. 

\begin{theorem}[ \cite{cluckers_non-archimedean_????}, Theorem 2.2.3]
\label{thm-global-lipsch}
Let $f :X \to K$ a definable function which is locally 1-Lipschitz. Then there is a $\kk$-definable partition of $X$ and a real $C\geq1$ such that on each part, $f$ is globally $C$-Lipschitz. 
\end{theorem}

\begin{theorem}
\label{thm-lipsch-dom-range}
Let $f : A\subseteq K^n \to K$ a definable function which is locally 1-Lipschitz. Then there is a finite definable partition of $A$ such that for each part $X$, the following holds. There is a coordinate projection $p : K^n \to K^{n-1}$ and some definable functions $g : X \to k^r$, $c : k^r\times K^{n-1} \to K$ and $d : k^r \times K^{n-1} \to K$ such that for each $\eta\in k^r$, the restrictions of $c(\eta,\cdot)$ and $d(\eta,\cdot)$ to $p(g^{-1}(\eta))$ are locally $1$-Lipschitz, and for each $w\in p(K^n)$, the set $g^{-1}(\eta)_w$ is in $c(\eta,w)$-config and the image of $g^{-1}(\eta)_w$ under $f_w$ is in $d(\eta,w)$-config. 
\end{theorem}

This theorem is not explicitly stated in \cite{cluckers_non-archimedean_????}, but it follows from the proof of Theorem 2.1.8 of \cite{cluckers_non-archimedean_????}. The idea is the same as Proposition \ref{partition-1lipsch} : use the Jacobian property on the centers of a cell decomposition to show that they are differentiable, then switch the order of the coordinates to bound the derivatives by 1, in the new coordinates, the centers are now locally 1-Lipschitz. 

In the residue characteristic zero case, Theorem 2.1.8 of \cite{cluckers_non-archimedean_????} shows that one can take $C=1$ in Theorem \ref{thm-global-lipsch}. In \cite{cluckers_approximations_2012}, Cluckers and Halupczok show that the same is true if $K$ is a finite extension of $\Qq_p$. We now generalize this result but only in the one dimensional case.

\begin{theorem}
\label{thm-global-lipsch-dim1}
Let $f :X\subseteq K^n\times Y \to K$ a definable function such that $X_y$ is one dimensional for every $y\in Y$ and the restriction of $f$ to $X_y$ is locally 1-Lipschitz. Then there is a $\kk$-definable partition of $X$  such that for any $y\in Y$, on each $\kk$-part $X_{y,\xi}$, $f$ is globally $1$-Lipschitz on $X_{y,\xi}$. 
\end{theorem}

To prove the theorem, we need the following lemma. 

\begin{lemma}
\label{zero-derivative}
Let $f : X\subseteq K \to K$ a definable differentiable function with $f'=0$. Then the image of $f$ is finite, there is a finite definable partition of $X$ in parts $X_i$ such that $f_{\vert X_i}$ is constant.
\end{lemma}
\begin{proof}
By the Jacobian property, $f$ is locally constant hence it follows from the definition of finite $b$-minimality.\end{proof}

If will also be convenient to use the following notion of thin cell.
\begin{definition}
Let $Y$ be definable, and $A\subseteq K\times Y$ be a 1-cell over $Y$. Then, for each $(t,y)\in A$, there is a unique maximal (for the inclusion) ball $B_{t,y}$ containing $t$ and satisfying $B_{t,y}\times \set{y}\subseteq A$.   For $y_0\in Y$, we call the collection of balls $\set{B_{t,y_0}}_{(t,y_0)\in A}$ the balls of the cell $A$ above $y_0$. 

We call  a 1-cell $A$ over $Y$ thin if the collection of balls of $A$ above any $y\in Y$ consists of at most one ball. 
\end{definition}

\begin{proof}[Proof of Theorem \ref{thm-global-lipsch-dim1}]
Suppose first that $n=1$. We have a definable function $f : X\subseteq K\times Y \to K$ such that $f(\cdot,y)$ is locally 1-Lipschitz on $X_y$ for every $y\in Y$. Apply Theorem \ref{thm-lipsch-dom-range} to $f$,  working on a part, call $c$ and $d$ the functions that parametrize the centers of the domain and range. After translating, assume $c=d=0$. Working on a $\kk$-part $X_{y,\eta}$, for some $\eta\in \kk$,  assume $X_{y,\eta}$ is a 1-cell (otherwise the result is trivial) with center 0, $f$ is injective on $X_{y,\eta}$, $f$ has the Jacobian property on $X_{y,\eta}$, and $f(X_{y,\eta})$ is a 1-cell with center 0. Fix $x_1,x_2\in X_{y,\eta}$. Note that if $\abs{x_1}=\abs{x_2}$, $i.e.$ they both lie in the same ball above $(y,\eta)$, then by the Jacobian property 
\[
\abs{\frac{\partial f}{\partial x}(x_1,y) \cdot (x_1-x_2)}=\abs{f(x_1,y)-f(x_2,y)},
\]
hence we are done because $\abs{\partial f(x_1,y)/\partial x}\leq 1$.

Assume now $\abs{x_1}\neq\abs{x_2}$. Then by our assumption, there exist $a_1,a_2, b_1, b_2\in \Zz$, $a_1\neq a_2$, $b_1\neq b_2$, $\ell, \ell'\in \Nn^*$, $\xi,\xi'\in \kk$ such that for $i=1,2$ the balls 
\[
B_{x_1,y}=\set{x\in K\mid \ord(x)=a_i, \ac_\ell(x)=\xi}
\]
are included in $X_{y,\eta}$, $x_i\in B_{x_i,y}$, and 
\[
f(B_{x_i,y})=\set{x\in K\mid \ord(x)=b_i, \ac_{\ell'}(x)=\xi'}.
\]
We have the equalities 
\[
\ord(f(x_1,y)-f(x_2,y))=\min_{i=1,2} b_i,
\]
\[
\ord(x_1-x_2)=\min_{i=1,2} a_i,
\]
and by the Jacobian property, comparing the sizes of the balls $B_{x_i,y}$ and $f(B_{x_i,y})$,
\[
\ell + \ord( \partial f(x_i,y)/\partial x)+a_i=\ell'+b_i.
\]
Hence we are done if $\ell'\leq \ell$, or if  $\ord(\partial f(x_i,y)/\partial x)\geq \ell'-\ell$. Therefore by taking a finite partition of $X$, we can assume  $\ord(\partial f(x,y)/\partial x)$ is constant on $X_{y,\eta}$. If $X_{y,\eta}$ is a thin cell, then once again we are done by the Jacobian property. By injectivity versus constancy (Proposition \ref{inj-cons}), up to taking another $\kk$-partition, we can assume $\partial f/\partial x$ is either constant or injective on $X_{y,\eta}$. 

Assume first that $\partial f /\partial x$ is constant on $X_{y,\eta}$. Then the function 
\[(x,y)\mapsto f(x,y)-\frac{\partial f(x,y)}{\partial x} \cdot x\]
has derivative 0 on $X_{y,\eta}$, hence it has a finite image by Lemma \ref{zero-derivative}, so by refining the partition, we can assume it is constant on $X_{y,\eta}$. Hence on $X_{y,\eta}$, 
\[
f(x,y)=a(y)x+b(y),
\]
for $a, b$ some definable functions. As $\abs{a(y)}\leq 1$, $f$ is then 1-Lipschitz on $X_{y,\eta}$.  

Consider now the last case, where $\partial f/\partial x$ is injective on $X_{y,\eta}$. By Theorem \ref{thm-lipsch-dom-range} applied to $\partial f/\partial x$ (for $n=1$ the Lipschitz assumption is superfluous), we get a refinement of the $\kk$-partition, and some definable functions $d,e$ such that each $\kk$-part $X_{y,\eta}$ is a 1-cell in $d(y)$-config, its image by $\partial f/\partial x$ is in $e(y)$-config, and $\partial f/\partial x$  has the Jacobian property on $X_{y,\eta}$. Once again, assume $X_{y,\eta}$ is not a thin cell otherwise we are done. By injectivity of $\partial f/\partial x$ , and the fact that $\ord(\partial f/\partial x)\geq 0$, up to removing finitely many thin cells that we can treat separately, we can assume $\ord(\partial f(x,y)/\partial x- e(y))\geq \alpha$, where $\alpha$ is an integer chosen as follows. Consider the function 
\[
g(x,y)=f(x,y)-e(y)\cdot x.
\]
Apply the beginning of the proof to the function $g$. We get a refinement of the decomposition, and call $\tilde \ell, \tilde \ell'$ the depth of embedding of the decomposition adapted to $g$. It follows that $g$ is globally 1-Lipschitz on the $\kk$-parts if $\ord(\partial g/\partial x)\geq \tilde \ell'-\tilde \ell$. Now set $\alpha=\tilde \ell'-\tilde \ell$, as 
\[
\frac{\partial g}{\partial x}(x,y)=\frac{\partial f}{\partial x}(x,y) -e(y),
\]
$g(\cdot,y)$ is globally 1-Lipschitz on $X_{y,\eta}$. The function $f$ is then a sum of two globally 1-Lipschitz functions of $X_{y,\eta}$ (because $\ord(e(y))\geq 0)$, hence by ultrametric inequality, it is itself globally 1-Lipschitz.

The general case reduces to the previous one, using Proposition \ref{partition-1lipsch} on $X$ and the chain rule for derivation. 
\end{proof}

\subsection{Curve selection}

We prove a weak form of the curve selection lemma. A full curve selection lemma is proved in the $p$-adic semi-algebraic case in \cite{scowcroft_structure_1988} and sub-analytic case in \cite{denef_p-adic_1988}. In our context, because the residue field can be infinite, there are no Skolem functions, hence definable choice is not available. Hence we cannot find a truly definable arc, however, we can find some uniformly definable family of arcs. Our result is weak in the sense that the arcs we get are Lipschitz continuous, but not analytic, this is the drawback of working with an abstract theory. In the the semi-algebraic language for equicharacteristic zero fields, Nowak shows in \cite{nowak_results_????} existence of definable analytic arcs, using resolution of singularities. 

Before stating the result, we prove some lemmas about limit values of Lipschitz definable functions needed for the proof. 

\begin{lemma}
\label{lemma-adh-value}
Let $f : U \to K$ a bounded function, $x\in \overline{U}$ and $\Gamma$ the graph of $f$. Then $\overline{\Gamma}\cap (\set{x_0}\times K^{n-d})$ is non-empty. 
\end{lemma}

\begin{proof}
One can assume $x=0$, and let $\beta$ a bound from below for $\ord(f)$. Let \[X=\set{(z,\gamma)\in K\times \Zz\mid \exists u\in U, \ord(u)\geq \gamma, \ord(f(u)-z)\geq \gamma}.\] Then $X$ is non-empty and projects surjectively on its second coordinate. By tameness, $X$ can be assumed of the form 
\[\set{(z,\gamma)\in K\times \Zz\mid \ac(z-c)\in D_1, \ord(z-c)\in D_2(\gamma)},\]
where $D_1$ is a definable of $k$, and $D_2(\gamma)$ is a $\gamma$-definable of $\Zz$. By definable choice for Presburger sets, there is a definable Presburger function $\alpha : \Zz \to \Zz$ with $\alpha(\gamma)\in D_2(\gamma)$ for every $\gamma\in \Zz$. 

Let show by contradiction that $\alpha(\gamma)$ is bounded from below when $\gamma$ goes to $+\infty$. If not, we can find a $\gamma$ with 
$\alpha(\gamma')<d< \gamma$,
for $d=\min\set{\ord(c),\beta}$.
By the description of $X$ above, there is a $z\in X$ with $\ord(z-c)=\alpha(\gamma)$. Then by definition of $X$, there is a $u\in U$ such that $\ord(u)\geq \gamma$, $\ord(f(u)-z)\geq \gamma$. Using the ultrametric inequality, we get 
\begin{eqnarray*}
\alpha(\gamma)=\ord(z-c)&=&\ord(z-f(u)+f(u)-c)\\
&\geq& \min\set{\ord(f(u)-z),\ord(f(u)-c)}\geq \min\set{\gamma,\beta}\geq d.
\end{eqnarray*}
Contradiction, hence $\alpha(\gamma)$ is bounded from below when $\gamma$ goes to $+\infty$. Then one can extract a sequence of $\gamma$ such that either $\alpha(\gamma)=\alpha_0$ is constant, either $\alpha(\gamma)$ goes to $+\infty$. In the first case, any $z$ with $\ord(z-c)=\alpha_0$ satisfy $(0,z)\in \overline{\Gamma}$, in the second one, $(0,c)\in \overline{\Gamma}$. 
\end{proof}

\begin{lemma}
\label{finite_acc_lemma}
Let $f: U\subseteq K^d\to K^{n-d}$ a globally $C$-Lipschitz function, of graph $\Gamma$. Then for any $x_0\in  \overline{U}$, there are only finitely many points of $\overline{\Gamma}$ which projects to $x_0$ under the coordinate projection $K^d\times K^{n-d}\to K^d$, $i.e.$ $X:=\overline{\Gamma}\cap \set{x_0}\times K^{n-d}$ is finite. 

\end{lemma}
For the proof we will use the tangent cone that will be studied in more details in section \ref{subsection_tangent_cone}. The $K^\times$-tangent cone of a definable subset $Y$ of $K^n$ is the set $C_x^{K^\times}(X)$ defined by
 \[
\set{u\in K^n\mid \forall i\in \Zz, \exists y\in X, \exists \lambda\in K^\times, \ord(y-x)\geq i, \ord(\lambda(y-x)-u)\geq i}.
\]

\begin{proof}
Suppose by contradiction that $\dim(X)\geq 1$. Then $X$ contains a one dimensional ball. For a $z$ in such a ball, $\dim(C_z^{K^\times}(X))\geq 1$.

As $f$ is globally $C$-Lipschitz, if $z=(z_0,z_1)$, then for any $x\in U$, $\abs{f(x)-z_1}\leq C\abs{x-z_0}$, therefore 
\[C_z^{K^\times}(\Gamma)\subseteq\set{(x,y)\in K^d \times K^{n-d} \mid \abs{y}\leq C\abs{x}}.\] 

As $C_z^{K^\times}(X)\subseteq C_z^{K^\times}(\overline{\Gamma})\cap C_z^{K^\times}(\set{x_0}\times K^{n-d})$, we then have
 
\[C_z^{K^\times}(X)\subseteq \set{(x,y)\in K^d \times K^{n-d} \mid \abs{y}\leq C\abs{x}} \cap (\set{0}\times K^{n-d})=\set{0},\]
 contradiction. 
\end{proof}

\begin{corollary}
\label{finite_accumulation}
Let $f: U\subseteq K^n\to K$ a locally $1$-Lipschitz function. Then $f$ has a nonzero finite number of limit values at any $z\in \overline{U}$. 
\end{corollary}

\begin{proof}
By Theorem \ref{thm-global-lipsch}, there is a $C\geq1$ and $\kk$-partition of $U$ such that $f$ is globally $C$-Lipschitz on each of the $\kk$-parts of $U$. For any $z\in \overline{U}$, by the following Lemma \ref{lem-kpart-closure} there is a $\kk$-part such that $z$ is in the closure of this $\kk$-part. Moreover, as any definable bijection between the residues sorts and the valued field sort is finite, it suffices to show that the result holds on each of the $\kk$-parts. But then $f$ is globally $C$-Lipschitz, hence bounded on some neighborhood of $z$, then by Lemma \ref{lemma-adh-value} $f$ has a limit value at $z$. By Lemma \ref{finite_acc_lemma}, $f$ has only a finite number of limit values. 
\end{proof}

\begin{lemma}
\label{lem-kpart-closure}
Let $X\subseteq K^n$ be definable and $g : X\to A$ a $\kk$-partition of $X$. Then the closure of $X$ is equal to the union of the closures of its $\kk$-parts.
\end{lemma}
\begin{proof}
Fix $u\in \overline{X}$. Consider the $u$-definable set
\[
D=\set{(n,\xi)\in \Zz\times A \mid \exists z\in X,\ord(z-u)\geq n, g(z)=\xi}.
\]
By definition of $u$, $D$ projects surjectively on $\Zz$. Because the theory splits, $D$ is a finite union of cartesian products of definable sets of $\Zz$ and $A$, hence there is a $\xi\in A$ such that for all $n\in \Zz, (n,\xi)\in D$, which means  $u\in \overline{g^{-1}(\xi)}$.
\end{proof}

\begin{proposition}[Weak curve selection]
\label{CSL}
Let $X\subseteq K^n$ be some definable, and $a\in \overline{X}\setminus X$. Then for some integer $e\geq1$,  there exists a definable function $\sigma : B_e\times D \to X\cup\set{a}$ for some nonempty definable $D$, where $B_e=\set{x\in K\mid \abs{x}\leq1, \ac_e(x)=1}$ and $B_e^*=B_e\setminus \set{0}$, and such that if we denote for any $d\in D$
\[
\sigma_d : t\in B_e \mapsto \sigma(t,d),
\]
the function $\sigma_d$ is 1-Lipschitz, $\sigma_d(0)=a$ and $\sigma_d(B_e^*)\subseteq X$. In the tame case, that is, when the residue characteristic is zero, one can take $e=1$.
\end{proposition}

To simplify the notations we will write down the proof in the equicharacteristic zero case, with $B=B_1$. For the mixed case, replace $k=k_K$ by some residue rings, and at the very end, $e$ is the greatest $s$ that appears in the $R_{s,K}^s$. 

\begin{proof}
Observe first that we only need to show the existence of a locally 1-Lipschitz function $\sigma$ satisfying the other requirements of the proposition, then use Theorem \ref{thm-global-lipsch-dim1} to get a globally 1-Lipschitz function. 

We proceed by induction over $n$. 
Up to translating, assume $a=0$. By Theorem \ref{thm-lipsch-dom-range}, there is a finite partition of $X$ into definable parts, such that on each part $A$, the following holds. There exists an $s\geq 0$, a coordinate projection $p : K^n \to K^{n-1}$, and some $\caL$-definable functions 
\[
g : A \to k^s \mrm{and} c : k^s \times K^{n-1} \to K
\]
such that for all $\eta \in k^s$, the restriction of $c(\eta, \cdot)$ to $(g^{-1}(\eta))$ is locally $1$-Lipschitz, and for each $w$ in $p(K^n)$, the set $g^{-1}(\eta)_w$ is in $c(\eta, w)$-config. There is a part $A$ such that $a\in \overline{A}$, and one can assume $X=A$. If $p(a)\in p(X)$ and $p(a)$ is isolated in $p(X)$, one is done by using the induction hypothesis with $n=1$ for $p^{-1}(a)\cap X$. Otherwise up to reducing $X$, one can assume that $p(a)\in \overline{p(X)} \setminus p(X)$. 

Let's fix an $\eta$ such that $g^{-1}(\eta)$ is nonempty, $a\in\overline{g^{-1}(\eta)}$ and $p(a)\in \overline{p(g^{-1}(\eta))} \setminus p(g^{-1}(\eta))$. Such an $\eta$ exists because the theory splits. Up to reducing $X$ again, using Corollary \ref{finite_accumulation} one can assume $\lim_{w\to 0} c(\eta,w)=c_0$ exists. One can assume $X=g^{-1}(\eta)$. The definable $X$ is now $\eta$-definable, and $\eta$ will be one of the parameters of the function we are defining.

By the induction hypothesis with $n-1$, there is a definable with parameters $1$-Lipschitz function $\sigma : B \to p(X)\cup\set{a}$ such that $\sigma(0)=p(a)$ and $\sigma(B\setminus \set{0})\subseteq p(X)$. 

If $c_0\neq 0$, as $0\in \overline{X}$, for all $w\in p(X)$ $0\in X_w$, hence the function $\sigma'=(\sigma,0)$ answer the problem. 

Let's consider the case $c_0= 0$.  We know that $X$ is in $c(\eta,w)$-config. If it is in graph-config, then  the function $\sigma'=(\sigma, c\circ\sigma)$ works. Otherwise 
\[X_w=\set{x\in K\mid \ac(x-c(\eta,w))=\eta_{i_0}, \ord(x-c(\eta,w))\in D_2},\]
 with $D_2$ an $(\eta,w)$-definable of $\Gamma$ and $\eta_{i_0}$ is one of the coordinates of $\eta$. As the theory splits, $D_2$ is in fact $w$-definable, denote it by $D_2(w)$. By Presburger quantifier elimination and compactness, up to reducing $X$ again, one can assume $D_2(w)$ is a finite union of cells of the form $\set{\alpha\in \Gamma\mid \alpha_1(w)\square_1 \alpha\square_2 \alpha_2(w), \alpha\equiv_e i}$, for some definable functions $\alpha_1,\alpha_2$,  integers $e,k$ and $\square_j$ either $\leq$ or and empty condition, and $\alpha_j(w)\in D_2(w)$. As $0\in\overline{g^{-1}(\eta)}$, there is a cell with $\alpha_2$ unbounded when $w$ goes to $0$ or with the condition $\square_2$ empty. 

Suppose we are in the first case. Let $\Delta$ the projection on the last two variables of 
\[\set{(w,\alpha,\beta)\in p(X)\times \Zz^2\mid \abs{w}=\alpha, \alpha_2(w)=\beta}.\] Because $\alpha_2(w)$ is unbounded when $w$ goes to 0, the point $(\infty,\infty)$ is in the closure of $\Delta$ for the order topology. By Presburger cell decomposition and existence of Skolem functions in Presburger arithmetic, there exists integers $\alpha_0,a_0, b_0, e_0,i_0$ with $a_0\geq1$ and $e_0\geq 1$ such that $\Delta$ contains the definable 
\[\set{(\alpha,\beta)\in\Zz^2\mid \alpha_0\leq \alpha, \alpha\equiv_{e_0} i_0, \beta=a_0\frac{\alpha-i_0}{e_0}+b_0}.\]
By the induction hypothesis there is a definable with parameters $1$-Lipschitz function $\sigma : B \to p(X)\cup\set{p(a)}$ such that $\sigma(0)=p(a)$ and $\sigma(B\setminus \set{0})\subseteq p(X)$. By tameness, the function $t\mapsto\ord(\sigma(t))$ is of the form $f(\ord(t))$, with $f$ a Presburger function, hence there are some $\alpha_1,a_1, b_1, e_1,i_1$ with $a_1\geq1$ and $e_1\geq 1$ such that for $\alpha_1\leq\ord(t)$ and $\ord(t)\equiv_{e_1}i_1$, $\ord(\sigma(t))=a_1\frac{\ord(t)-i_1}{e_1}+b_1$, and by composition, there are some integers $i_2,e_2, a_2, b_2$ with $a_2\geq 1$, $e_2\geq1$ such that for all $t\in B$ and $\sigma'(t):=\sigma(\varpi_K^{i_2}t^{e_2})$, $\alpha_2(\sigma'(t))=a_2\ord(t)+b_2$. Choose a $d\in X_{\sigma'(1)}$. Then the function $t\mapsto (\sigma'(t), c(\eta,\sigma'(t))+dt^{a_2})$ maps 0 to 0 and $B\setminus \set{0}$ to $X$ and is locally $1$-Lipschitz. The latter case, where the cell of $D_2(w)$ is of the form $\alpha_1(w)\leq \alpha$ is similar.
\end{proof}

By working over parameters and using logical compactness, we get the following parametrized curve selection lemma. 
\begin{proposition}[Parametrized curve selection]
\label{CSL-par}
Let $X\subseteq K^n\times Y$ be definable and for every $y\in Y$,  a $y$-definable point $a_y\in \overline{X_y}\setminus X_y$. Then for some integer $e\geq1$,  there exists a definable function $\sigma : B_e\times D\times Y \to K^n$ for some nonempty definable $D$, where $B_e=\set{x\in K\mid \abs{x}\leq1, \ac_e(x)=1}$ and $B_e^*=B_e\setminus \set{0}$, and such that if we denote for any $(d,y)\in D\times Y$
\[
\sigma_{d,y} : t\in B_e \mapsto \sigma(t,d,y),
\]
the function $\sigma_d$ is 1-Lipschitz, $\sigma_d(0)=a_y$ and $\sigma_d(B_e^*)\subseteq X_y$. In the tame case, that is, when the residue characteristic is zero, one can take $e=1$.
\end{proposition}

The following immediate corollary will be useful in the study of regular stratifications. In addition to the use of the curve selection, one has to notice that any definable function is differentiable almost everywhere, that is, everywhere but on a definable set of lower dimension, which is guaranteed by the Jacobian property. 
\begin{corollary}
\label{cor-CSL}
Let $U\subseteq K^m$ be a nonempty open definable and $S\subseteq  K^n$, such that for all $u\in U$, $(u,0)\in \overline{S_u}\setminus S_u$. 
Then there is a nonempty open $V\subseteq U$ $\kk$-definable and a definable with parameters 
\begin{eqnarray*}
\overline{\sigma} : V\times B_e &\to& Y\cup (V\times \set{0}^{n-m})\\
(u,s)&\mapsto& (u,\sigma(u,s))
\end{eqnarray*}
such that $\overline{\sigma}(u,0)=0$, $\overline{\sigma}(V\times B_e^*)\subseteq V$,
for any fixed $u\in V$, $\sigma(u,\cdot)$ is $1$-Lipschitz, and for any $s\in B_e$, $\sigma(\cdot,s)$ is differentiable. 
\end{corollary}

The following proposition will be of crucial use in the study of tangent cones. In the $p$-adic case, the analogous result is proven using analytic curve selection and Whitney's lemma. As we do not have analytic curve selection lemma, we cannot use Whitney's lemma but we will instead use our Theorem \ref{thm-global-lipsch-dim1} to get a globally Lipschitz function. 
\begin{proposition}
\label{prop-tangent-incl}
Let $f : X\subseteq K^d \to K^{n-d}$ a definable function,  $\Gamma$ be its graph and suppose $f$ is locally $C$-Lipschitz. Fix $\Lambda\in \Dd$ and $z=(z_1,z_2)\in \overline{\Gamma}$. Then there is an open neighborhood  $U$ of $z_1$ in $K^d$ such that the following inclusions hold :
\[
\set{(x-z_1,f(x)-z_2)\mid x\in U}\subseteq\set{(x,y)\in K^d \times K^{n-d} \mid \abs{y}\leq C\abs{x}}
\]
and
\[
C_z^{\Lambda}(\Gamma)\subseteq\set{(x,y)\in K^d \times K^{n-d} \mid \abs{y}\leq C\abs{x}}.
\]
\end{proposition}
For the definition of the tangent cone $C_z^{\Lambda}(\Gamma)$, see Definition \ref{def-tan-con}.

\begin{proof}
After translating and rescaling, one can assume $z=0$ and $C=1$ and it suffices to show the existence of a neighborhood $U$ of $z_1$ such that the first inclusion holds. Assume by contradiction that there is a  $\delta>0$ such that the set 
\[
\Gamma':=\set{(x,f(x))\mid x\in X, \abs{f(x)}\geq (1+\delta)\abs{x}}
\]
is nonempty and satisfies $0\in \overline{\Gamma'}\backslash \Gamma'$. By the curve selection lemma \ref{CSL}, there is a definable function $\sigma : B_e\times D \to K^n$ with $\sigma_d$ 1-Lipschitz, $\sigma_d(0)=0$ and $\sigma_d(B_e^*)\subseteq \Gamma'$. Denote $\tilde{\sigma}:B_e\times D \to K^d$ the first $d$ components of $\sigma$ and set
\[
X':=\set{(x,d)\in K^d\times D\mid \exists t\in B_e^*, x=\tilde{\sigma}(t,d)}.
\]
For any $d\in D$, the set $X'_d$ is one dimensional and $(0,0)$ lies in the closure of the graph of $f$ on $X'_d$. By Proposition \ref{thm-global-lipsch-dim1}, up to taking a $\kk$-partition, one can assume $f$ is globally 1-Lipschitz on $X'_d$. Fix a $x\in X'_d$. There is a $y\in X'_d$ with $\abs{y}<\abs{x}$ and $\abs{f(y)}\leq \abs{f(x)}$. Then by 1-Lipschitz continuity of $f$, 
\[
\abs{f(x)-f(y)}\leq \abs{x-y},
\]
but $\abs{x}=\abs{x-y}$ and $\abs{f(x)}=\abs{f(x)-f(y)}$, hence $\abs{f(x)}\leq\abs{x}$. Contradiction with $x\in \Gamma'$. 
\end{proof}

The following proposition is a partial converse to Proposition \ref{prop-tangent-incl}. 

\begin{proposition}
\label{prop-tangent-lipschitz}
Let $f : U \to K^{n-d}$ a differentiable definable function on a nonempty open $U\subseteq K^d$, and let $\Gamma$ be its graph. Assume 

\[C_{(x_0,f(x_0))}^{\Lambda}(\Gamma)\subseteq\set{(x,y)\in K^d \times K^{n-d} \mid \abs{y}\leq C\abs{x}}.
\]
Then for all $i=1,...,d$, $\abs{\partial f/\partial x_i (x_0)}\leq C$, and in particular, $f$ is locally $C$-Lipschitz. 
\end{proposition}
\begin{proof}
For $i=1,...,d$, the vector $(0,...,0,1,0,..., \partial f/\partial x_i (x_0))$ (with a 1 in $i$-th position) is in $C_{(x_0,f(x_0))}^{\Lambda}(\Gamma)$. 
\end{proof}

\subsection{Recap on motivic integration}

We use Cluckers-Loeser theory of motivic integration defined in \cite{cluckers_motivic_????}, see also \cite{cluckers_constructible_2008} for the residue characteristic zero case. We review here the main steps of the construction in order to fix the notation. 

Fix either a tame theory or a mixed tame theory $\mathcal{T}$ in a language $\mathcal{L}$ and $X$ a definable subassignment.

\paragraph{Value group}
Let $\Aa$ be the ring
\[ 
\Aa:= \Zz\left[\eL,\eL^{-1},{\left({\frac{1}{1-\eL^{-n}}}\right)_{n\geq 1}}\right]. 
\]
For every real $q>1$, define $\theta_q : \Aa \to \Rr$ to be the ring morphism generated by $\eL \mapsto q$, and $\Aa$ the sub-semiring $\Aa_+ :=\set{a \in \Aa \mid  \forall q>1, \theta_q(a)\geq 0}$. The ring $\Aa$ is equipped with a partial order defined by $a\leq b$ if and only if $a-b \in \Aa_+$. 

The application degree $\Zz[\eL] \to \Nn$ extends uniquely as a map $\deg_\eL : \Aa \to \Zz$ satisfying $\deg_\eL(ab)=\deg_\eL(a)+\deg_\eL(b)$. A sequence $(a_n)_n\geq 0$ of elements of $\Aa$ is said to converge to 0 if $\deg_\eL(a_n)$ goes to $-\infty$. This defines a topology on $\Aa$. 

Let $\Pp(X)$ the ring generated by functions 
\begin{itemize}
\item $\alpha : X \to \Zz$ definable,
\item $\eL^\alpha$, with $\alpha : X \to \Zz$ definable,
\item constants of $\Aa$.
\end{itemize}
The sub-semiring $\Pp_+(X)$ of $\Pp(X)$ is constituted of functions with value in $\Aa_+$. 

In the definition above, we allow a slight abuse of notation by identifying $\Zz$ with the value group sort, but it is harmless since we only consider models of $\mathcal{T}$ that are $\mathcal{T}$-fields, hence with value group $\Zz$. 

The partial order on $\Aa$ induces a partial order on $\Pp_+(X)$ and $\Pp(X)$ :
for $\varphi,\psi\in \Pp(X)$, define $\varphi\geq \psi$ if $\varphi-\psi\in \Pp_+(X)$.

For any $\mathcal{T}$-field $K$, $\varphi\in \Pp(X)$, define $\theta_{q,K}(\varphi)$ to be the function
\[
\theta_{q,K}(\varphi) : x\in X(K) \mapsto \theta_q(\varphi(x)).
\]
A function $\varphi\in \Pp(X[0,0,s])$ is said to be $X$-integrable if for any $\mathcal{T}$-field $K$ and any $x\in X(K)$ the family $(\theta_{q,K}(\varphi)(x,i))_{i\in \Zz^s}$ is absolutely summable (as a family of real numbers). 

Denote by $I_X\Pp(X[0,0,s])$ (resp. by $I_X\Pp_+(X[0,0,s])$) the $\Pp(X)$-module (resp. the $\Pp_+(X)$-semimodule) of $X$-integrable $\Pp(X)$ functions (resp. $\Pp_+(X)$ functions). 

\begin{proposition}[\cite{cluckers_motivic_????}, Theorem 5.1]
Let $\varphi\in \Pp(X[0,0,s])$ an $X$-integrable function, then there is a unique $\mu_{/X}(\varphi)\in \Pp(X)$ such that for all $q>1$, $K$ a $\mathcal{T}$-field and $x\in X(K)$, 
\[
\theta_{q,K}(\mu_{/X}(\varphi))(x)=\sum_{i\in \Zz^s} \theta_{q,K}(\varphi)(x,i).
\]
Moreover, the map $\varphi \mapsto \mu_{/X}(\varphi)$ is a morphism of $\Pp(X)$-modules
\[
\mu_{/X} : I_X \Pp(X[0,0,s])\to \Pp(X)
\]
and send $I_X \Pp_+(X[0,0,s])$ to $\Pp_+(X)$. 
\end{proposition}

Define also $\Pp^0(X)$ (resp. $\Pp_+^0(X)$) as the sub-ring of $\Pp(X)$ (resp. as the sub-semiring of $\Pp_+(X)$) generated by characteristic functions $\mathbf{1}_Y$ of every definable $Y\subseteq X$ and the constant function $\eL-1$. 

For any $f : X\to Y$ definable, define the pullback by
\begin{eqnarray*}
f^* : \Pp_+(Y)&\to& \Pp_+(X)\\
\varphi &\mapsto& \varphi\circ f.
\end{eqnarray*}

\paragraph{Residue ring}
Let $\Qm_+(X)$ the quotient of the free abelian semigroup generated by symbols $[Y_{/X}]$, where $Y\subseteq X[0,r,0]$ is definable for some $r=(r_1,...,r_s)\in \Nn^s$, by the relations :
\begin{itemize}
\item $[Y_{/X}]=[Y'_{/X}]$ if $Y$ and $Y'$ are definably isomorphic over the projection to $X$,
\item $[Y_{/X}]+[Y'_{/X}]=[(Y\cup Y')_{/X}]+[(Y\cap Y')_{/X}]$,
\item $[\emptyset_{/X}]=0$
\item $[Y_{/X}]=[Y'_{/X}]$ if for some definable $Z$ of $X[0,r,0]$, one has 
\[Y'=p^{-1}(Z) \mrm{and} Y=Z[0,1,0],\]
where $p : X[0,(r_1+1,r_2,...,r_s),0]\to X[0,r,0]$ is the natural projection between residue rings. 
\end{itemize}
In particular, the classes of $X[0,r,0]$ and $X[0,r',0]$ in $\Qm_+(X)$ are identified for any tuples $(r_i)_i$ and $(r'_j)_j$ such that $\sum_i r_i=\sum_j r'_j$. Note also that the last relation is vacuous if the residue field is perfect. 

The semigroup $\Qm_+(X)$  carries a semiring structure with multiplication for subsets $Y\subseteq X[0,m,0]$ and $Y'\subseteq X[0,m',0]$ given by
\[
[Y_{/X}]\cdot [Y'_{/X}]:=[Y\times_XY'_{/X}].
\]

We denote $\Qm(X)$ the groupification of $\Qm_+(X)$. Then the canonical morphism $\Qm_+(X)\to \Qm(X)$ has no reason in general to be injective. If the theory of the residue field is the theory of algebraically closed fields of characteristic zero, Karzhemanov shows in \cite{karzhemanov} that it is not injective.

To keep information at the level of residue rings, we define formally integration of residue variables. Let $r$ be a tuple of integers, the integral of $[Y_{/X}]\in \Qm_+(X[0,r,0])$ with respect to $X$ is 
\[
\mu_{/X} : [Y_{X[0,r,0]}\in \Qm_+(X[0,r,0]) \mapsto [Y_{/X}]\in \Qm_+(X).
\]

\paragraph{Motivic constructible functions}
We identify $[X\times \mathbb{G}_m(k)_{/X}]$ with $\eL-1$, and for $Y\subseteq X$ definable, $[Y_{/X}]$ with $\mathbf{1}_Y$, hence $\eL$ is identified with $[X\times\Aa_k^1]=[X[0,(1),0]]$. Therefore there are inclusions 
\[\Pp_+^0(X)\hookrightarrow \Qm_+(X) \mrm{and} \Pp^0(X)\hookrightarrow \Qm(X).\]

Define the semiring of positive motivic constructible functions as the tensor product
\[ \Cm_+(X):=\Pp_+(X)\otimes_{\Pp_+^0(X)} \Qm_+(X)\]
and the ring of motivic constructible functions
\[ \Cm(X):=\Pp_+(X)\otimes_{\Pp_+^0(X)} \Qm_+(X).\]

For any $f :X\to Y$ definable, the pullback morphisms on $\Pp_+(Y)$ and $\Qm_+(Y)$ extend canonically to a pullback morphism $f^* : \Cm_+(Y)\to \Cm_+(X)$. 

There is also a canonical morphism $ \iota : \Cm_+(X)\to  \Cm(X)$, but it has no reason to be injective.

Using the fact that the theory splits, we get canonical isomorphisms
\[ \Cm_+(X[0,r,s])\simeq \Pp_+(X[0,0,s])\otimes_{\Pp_+^0(X)}\Qm_+(X[0,r,0]) \]
and 
\[ \Cm(X[0,r,s])\simeq \Pp(X[0,0,s])\otimes_{\Pp^0(X)}\Qm(X[0,r,0]). \]

These isomorphisms allow one to define canonically integrals over auxiliary variables. Define 
\[I_X\Cm_+(X[0,r,s]):= I_X\Pp_+(X[0,0,s])\otimes_{\Pp_+^0(X)}\Qm_+(X[0,r,0])\]
and 
\begin{eqnarray*}
\mu_{/X} : I_X\Cm_+(X[0,r,s])&\to&\Cm_+(X)\\
a\otimes b &\mapsto& \mu_{/X}(a)\otimes\mu_{/X}(b).
\end{eqnarray*}

To integrate over valued field variables, one needs to work component by component, using cell decomposition. 
Fix $K$ a $\mathcal{T}$-field, for any ball $B=a+b\mathcal{O}_K\subseteq K$ and real $q>1$, define $\theta_{q,K}(B)$ to be the real $q^{-\ord(b)}$. 
A countable disjoint union of balls of $K$, all of different radii, is called a step domain. A function $\varphi : K \to \Rr_{\geq 0}$ is called a step function if there is a unique step domain $S$ such that $\varphi$ is constant and non-zero on each ball of $S$ and zero outside $S\cup\set{a}$ for some $a\in K$. 

Let $q>1$ be a real. A step-function $\varphi : K \to \Rr_{\geq 0}$ on step-domain $S$ is said to be $q$-integrable if 
\[
\sum_{B\subseteq S}\theta_{q,K}\varphi(B)\leq +\infty,
\]
the summation being on the balls $B$ of $S$. This value is the $q$-integral of $\varphi$. 

Let $\varphi\in \Pp_+(X[1,0,0])$, we say $\varphi$ is an $X$-integrable step function if for any $\mathcal{T}$-field $K$, any $x\in X(K)$ and $q>1$, the function
\begin{eqnarray*}
\theta_{q,K}(\varphi)(x,\cdot) : K &\to& \Rr_{\geq 0}
t \mapsto \theta_{q,K}(\varphi)(x,t)
\end{eqnarray*}
is a $q$-integrable step-function. For such a function, there is a unique function $\psi\in \Pp_+(X)$ such that for any $x\in X(K)$, $\theta_{q,K}(\psi)(x)$ equals the $q$-integral of $\theta_{q,K}(\varphi)(x,\cdot)$.  We say that $\varphi$ is $X$-integrable and set
\[
\mu_{/X}(\varphi):=\psi,
\]
the integral of $\varphi$ over $X$. 

We can now define the integral of a general motivic constructible function in relative dimension 1. 

\begin{lemma}[\cite{cluckers_motivic_????}, Lemma-Definition 8.2]
Let $\varphi\in\Pp_+(X[1,0,0])$. Say that $\varphi$ is $X$-integrable if there is a $\psi\in\Pp_+(X[1,r,0])$ with $\mu_{/X[1,0,0]}(\psi)=\varphi$ such that $\psi$ is $X[0,r,0]$-integrable in the previous sense and set
\[
\mu_{/X}(\varphi)=\mu_{/X}(\mu_{/X[0,r,0]}(\psi))\in \Cm_+(X)
\]
This definition is independent of the choices, it is the integral of $\varphi$ over $X$. 
\end{lemma}
Using cell decomposition, for any $\varphi\in \Cm_+(X[1,0,0])$ one can find a step function $\psi \in \Pp_+(X[1,r,0])$ for some tuple $r$ such that $\mu_{/X[1,0,0]}(\psi)=\varphi$, but one needs to show that this definition is independent of the choice of $\psi$. 

One can now define general integration by induction on the relative dimension. 

\begin{lemma}[\cite{cluckers_motivic_????}, Lemma-Definition 9.1]
Let $\varphi\in \Cm_+(Z)$, for some $Z=X[n,r,s]$ and $X$ definable. Say $\varphi$ is $X$-integrable if there is a definable $Z'\subseteq Z$ of complement in $Z$ of relative dimension $<n$ over $X$ and a choice of the order of coordinates on $Z$ such that $\varphi':={\mathbf 1}_{Z'}\cdot\varphi$ is $X[n-1,r,s]$-integrable and $\mu_{/X[n-1,r,s]}(\varphi')$ is $X$-integrable. In this case, set
\[
\mu_{/X}(\varphi):= \mu_{/X}(\mu_{/X[n-1,r,s]}(\varphi'))\in \Cm_+(X).
\]
This definition is independent of the choice of the order of coordinates. When $X=h[0,0,0]$ is the final object, one denote $\mu:=\mu_{/X}$. 
\end{lemma}

We can now integrate functions  $\varphi\in\Cm_+(h[n,r,s])$, but the integral $\mu(\varphi)$ is zero if (the support) of $\varphi$ is of dimension $<n$. By analogy with classical measure theory, one now define the integral of dimension $n$ of functions embedded in higher dimensional spaces, using a change of variable formula.  
Fix $K$ a $\mathcal{T}$-field. For any $f : A\subset K^n \to K^n$, let $\Jac f : A \to K$ the determinant of the Jacobian matrix when it is defined, 0 otherwise. In the relative case, for a function $f : A\subset C\times K^n \to C\times K^n$ such that the diagram of projections to $C$ commutes, one define $\Jac_{/C} : A \to K$ to be the function  satisfying for every $c\in C$, 
\[
(\Jac_{/C} f)(c,z)=(\Jac f_c)(z),
\]
with $f_c$ the function sending $z$ to $t$ such that $f(c,z)=(c,t)$. 

If $f$ is definable, the (relative) Jacobian is definable, by definability of partial derivatives, and one sets $\Jac_{/C} : A \to h[1,0,0]$ the unique definable function compatible for any $\mathcal{T}$-field $K$ with the definition above. 

One can now state the change of variable formula. The proof relies on the Jacobian property in dimension $1$. 

\begin{proposition}[\cite{cluckers_motivic_????}, Theorem 10.3]
Let $f : Z\subseteq X[n,0,0] \to Z'\subseteq X[n,0,0]$ a definable isomorphism over $X$, and $\varphi\in \Cm_+(Z')$. Then there is a definable $Y\subseteq Z$ with complement in $Z$ of dimension $<n$ relatively to $X$ and such that the relative Jacobian $\Jac_{/X} f$ is nonzero on $Y$. Then $\eL^{-\ord(\Jac_{/X} f)}f^*(\varphi)$ is $X$-integrable if and only if $\varphi$ is $X'$-integrable, and then
\[
\mu_{/X}(\eL^{-\ord(\Jac_{/X} f)}f^*(\varphi))=\mu_{/X}(\varphi)\in \Cm_+(X),
\]
with the convention $\eL^{-\ord(0)}=0$. 
\end{proposition}

We define now the graded ring of positive motivic constructible $\mathcal{F}$unctions. For any integer $d\geq 0$, define $\Cm_+^{\leq d}(X)$ to be the ideal of $\Cm_+(X)$ generated by characteristic functions $\mathbf{1}_Y$ for $Y\subseteq X$ of dimension $\leq d$. One also sets $\Cm_+^{\leq 1}(X)=\set{0}$. For any $d\geq 0$, set
\[\mathrm{C}_+^{d}(X):=\Cm_+^{\leq d}(X) /\Cm_+^{\leq d-1}(X),\]
\[\mathrm{C}_+(X):=\bigoplus_{d\geq 0} \mathrm{C}_+^{d}(X).\]
The nonzero elements of $\mathrm{C}_+^{d}(X)$ are functions of support of dimension $d$, defined up to sets of dimension $<d$. 

One can define integrals of $\mathcal{F}$unctions in $\mathrm{C}_+(X)$ that take into account dimension. Rather than a formal definition, we explain how to compute them. Suppose $X=h[n,0,0]$. If $[\varphi] \in \mathrm{C}_+^d(X)$. Using functoriality properties and working with parameters, one can suppose using Proposition \ref{partition-1lipsch} that $X$ is a graph of a 1-Lipschitz function $f : U\subseteq h[d,0,0] \to h[n-d,0,0]$, where $U$ is open in $K^d$. By abuse of notation, write $f : u\in U \mapsto (u,f(u))\in X$ and define
\[
\mu_d(\varphi):=\mu(f^*(\varphi)).
\]
This does not depend on the chosen representative of $[\varphi]$, because if $\psi\in \Cm_+(U)$ has support of dimension $<d$, then $\mu(\psi)=0$.

\section{Local density}
\label{section_local_density}

By analogy with complex and real cases, it would be tempting to define the local density of a definable $X\subseteq K^m$ at $x\in K^m$ by
\[
\lim_{n\to +\infty} \frac{\mu_d(X\cap S(x,n))}{\mu_d(S_d(0,n))},
\]
where $S(x,n)=\set{y\in K^m\mid \min_i\set{\ord(y_i-x_i)}=n}$ is the (non-archimedean) sphere of valuative radius $n$. 
However this limit does not exist in general, for example, if we take 
\[
X=\set{x\in K\mid \ord(x)=0 \mod 2}, 
\]
then $X\cap S(0,2k+1)=\emptyset$ and $X\cap S(0,2k)=S(0,2k)$, hence the sequence oscillate between 0 and 1. In this case, we will define the local density of $X$ to be 1/2. More generally, we will show that for every sequence of the form
\[
\theta_n= \frac{\mu_d(X\cap S(x,n))}{\mu_d(S_d(0,n))},
\]
there is an integer $e$ such that the subsequences $(\theta_{ke+i})_{k\geq 0}$ converge to some $d_i$ and one can define the density to be the arithmetic mean $\frac{1}{e} \sum_{i=1}^e d_i$. Note that this phenomenon is not anecdotic, if we compute the density of the cusp $x^2=y^3$ at 0, the same oscillation appears, and we get a density of 1/2.

\subsection{Mean value at infinity}

Let $X, Y$ be definable. If $x\in X$, we write $i_x$ the inclusion $\set{x}\times Y \to X\times Y$. Recall that there is a canonical pullback morphism  $i_x^* : \Pp_+(X\times Y)\to \Pp(Y)$, we will use sometimes the convenient notation $\varphi(x)=i_x^*(\varphi)$.

We define also $\widetilde{\Pp_+}(X):=\Pp_+(X)\otimes \Qq$ and $\widetilde{\Cm_+}(X):=\Cm_+(X)\otimes \Qq$.

\begin{definition}
We say that a function $\varphi\in \Pp_+(\Zz)$ admits a mean value at infinity, denoted by $\MV(\varphi)$ if there is a positive integer $e$ such that for all $i=0,...,e-1$, the subsequence $\varphi(i+ne)_{n\in \Nn}$ converges to a $d_i\in \Aa_+$. One sets $\MV(\varphi)=\frac{1}{e}\sum_{i=0}^{e-1}d_i\in \Qq\otimes_\Zz \Aa_+$. 

More generally, we say that a function $\varphi\in \Pp_+(X\times \Zz)$ has a mean value at infinity $\MV(\varphi)\in \widetilde{\Pp_+}(X)$ if and only if for every $x\in X$, $i_x^*(\varphi)$ has a mean value at infinity and one set $\MV(\varphi)\in \Pp_+(X)$ such that for every $x\in X$, $i_x^*(\MV(\varphi))=\MV(i_x^*(\varphi))$. 
\end{definition}

\begin{definition}
Let $X$ be definable. A function $\varphi\in\Pp_+(\Zz)$ is said to be bounded if there is an $M\in \Aa_+$ such that $\varphi<M$. A function $\varphi\in\Pp_+(X\times\Zz)$ is said to be $X$-bounded if for all $x\in X$, $i_x^*(\varphi)$ is bounded. 
\end{definition}

\begin{lemma}
\label{lem-MVPP}
Let $\varphi\in \Pp(X\times \Zz)$ an $X$-bounded function. Then $\varphi$ admits a mean value at infinity. 
\end{lemma}

\begin{proof}
By Presburger cell decomposition, one can assume there are some $e,k\in \Nn^*$, a definable function $b : X \to \Nn$, $h_{i,j}\in \Pp_+(X)$, $a_i\in \Nn, b_i\in \Zz$  for $i=1,...,k$ and $j=0,...,e-1$ such that for all $x\in X$, $j\in \set{0,...,e-1}$, $n\equiv j \mod e$ and $n\geq b(x)$, one has
\[
i_{x,n}^*(\varphi)=\sum_{i=1}^k h_{i,j}(x)\left(\frac{n-j}{e}\right)^{a_i}\eL^{b_i(\frac{n-j}{e})}.
\]
One can moreover assume that the pairs $(a_i,b_i)$ are pairwise distinct, and that the functions $h_{i,j}$ are not identically zero (unless perhaps the one corresponding to (0,0)). As $i_x^*(\varphi)$ is bounded, for all $i$, $b_i\leq 0$, and if $b_i=0$, then $a_i=0$. Let $i_0$ such that $(a_{i_0},b_{i_0})=(0,0)$. One then has $\lim_{n\to +\infty} i_{x,j+ne}^*(\varphi)=h_{i_0,j}(x)$. This implies that $\varphi$ has a mean value at infinity which is 
\[
\MV(\varphi)=\frac{1}{e}\sum_{j=0}^{e-1}h_{i_0,j}\in \widetilde{\Pp_+}(X).
\]
\end{proof}

\begin{definition}
Let $\varphi\in \Cm_+(X\times \Zz)$. One says that $\varphi$ has a mean value at infinity if it can be written $\varphi=\sum_i a_i\otimes b_i$ with $a_i\in \Qm_+(X)$, $b_i\in \Pp_+(X\times \Zz)$, such that the $b_i$ have mean values at infinity. One sets 
\[
\MV(\varphi)=\sum_i a_i\otimes \MV(b_i)\in \widetilde{\Cm_+}(X).
\]
By using some refined partition of $X$, one can check that this definition does not depend on the chosen representation of $\varphi$. 

One says that $\varphi\in \Cm_+(X\times \Zz)$ is $X$-bounded if it can be written $\varphi=\sum_i a_i\otimes b_i$ with $a_i\in \Qm_+(X)$, $b_i\in \Pp_+(X\times \Zz)$, such that the $b_i$ are $X$-bounded.
\end{definition}

\begin{lemma}
\label{lem-MVCM}
Let $\varphi\in \Cm_+(X\times \Zz)$ an $X$-bounded function. Then $\varphi$ admits a mean value at infinity. 
\end{lemma}

\begin{proof}
It follows from Lemma $\ref{lem-MVPP}$ and the definitions, using the canonical isomorphism
\[
\Cm_+(X\times \Zz)\simeq \Qm_+(X)\otimes_{\Pp_+^0(X)}\Pp_+(X\times \Zz).
\]
\end{proof}

Let $A$ be definable, and $\varphi\in \Cm_+(A[m,0,0])$, locally bounded, of support of dimension $\leq d$ relatively to $A$. Let $p : A[2m,0,1] \to A[m,0,0]$ the projection, and 
\[\mathcal{S}=A\times \set{(x,y,z)\in h[m+m,0,1]\mid \min_i\set{\ord(y_i-x_i)}=n}.
\]
Then $\mathbf{1}_\mathcal{S} \cdot p^*(\varphi)\in \Cm_+(A[2m,0,1])$, and define
\[
\gamma_d(\varphi):=\mu_{d,A[m,0,1]}(\mathbf{1}_\mathcal{S}\cdot p^*(\varphi))\in \Cm_+(A[m,0,1]).
\]
Then one normalizes by the volume $(1-\eL^{-d})\eL^{-nd}$ of the sphere $S^d(0,n)$ of dimension $d$ and radius $n$. One defines $\theta_d(\varphi)\in \Cm_+(A[m,0,1])$ such that
\[
\theta_d(\varphi)(n)=\frac{\gamma_d(\varphi)(n)}{(1-\eL^{-d})\eL^{-nd}}.
\]

\begin{proposition-definition}
\label{def-MLD}
If $\varphi\in \Cm_+(A[m,0,0])$ is locally bounded, of support of dimension $\leq d$ relatively to $A$, then $\theta_d(\varphi)$ is well defined and admits a mean value at infinity. One sets 
\[
\Theta_d(\varphi):=\MV(\theta_d(\varphi))\in \widetilde{\Cm_+}(A[m,0,0]),
\]
the motivic local density of $\varphi$. 
\end{proposition-definition}

If $\varphi=\mathbf{1}_X$, and $x\in K^m$, one denote $\Theta_d(X,x):=i_x^*(\Theta_d(\mathbf{1}_X))\in  \widetilde{\Cm_+}(\set{x})$. 

\begin{proof}
As $\varphi$ is locally bounded, at least for $n$ big enough, $\mathbf{1}_\mathcal{S}p^*(\varphi)$ is $A[m,0,1]$-integrable, and $\theta_d(\varphi)$ is $A[m,0,0]$ bounded, hence by Lemma $\ref{lem-MVCM}$ it has a mean value at infinity. 
\end{proof}

In more explicit terms, one has simply
\[\Theta_d(X,x)=\MV\left(n\mapsto\frac{\mu_d(X\cap S^m(x,n))}{\mu_d(S^d(0,n))}\right).\]

The following lemma shows that in the definition of local density we can replace the spheres by balls. 

\begin{lemma}
We keep the notations of \ref{def-MLD}. Let 
\[
\mathcal{B}:=A\times \set{(x,y,z)\in h[m+m,0,1]\mid \min_i\set{\ord(y_i-x_i)}\geq n},
\] 
$\gamma'_d(\varphi):=\mu_{d,A[m,0,1]}(\mathbf{1}_\mathcal{B}\cdot p^*(\varphi))$, and $\theta'_d(\varphi)(n):=\frac{\gamma'_d(\varphi)(n)}{\eL^{-nd}}$. Then $\theta'_d(\varphi)$ has a mean value at infinity and 
\[
\MV(\theta'_d(\varphi))=\Theta_d(\varphi).
\]

\begin{proof}
The limit $\MV(\theta'_d(\varphi))$ exists for the very same reasons that $\MV(\theta_d(\varphi))$. One has for all $n\in \Zz$
\[
\gamma_d(\varphi)(n)=\gamma'_d(\varphi)(n)-\gamma'_d(\varphi)(n+1), 
\]
and 
\[\theta_d(\varphi)(n)=\frac{1}{1-\eL^{-d}}(\theta'_d(\varphi)(n)-\eL^{-d}\theta'_d(\varphi)(n+1)).\]
By linearity of $\MV$, one has the result. 
\end{proof}
\end{lemma}

Here are some basic properties of the local density. 
\begin{proposition} Let $\varphi, \psi \in \Cm_+(A[m,0,0])$ be locally bounded, of supports of dimension $\leq d$ relatively to $A$. Let $X,Y\subseteq K^m$ definable of dimension $d$, with $X\cap Y$ of dimension $<d$. 
\begin{itemize}
\item $\Theta_d(\varphi+\psi)=\Theta_d(\varphi)+\Theta_d(\psi)$,
\item $\Theta_d(X\cup Y)=\Theta_d(X)+\Theta_d(Y)$,
\item $\Theta_d(X)=\Theta_d(\overline{X})$.
\end{itemize}
\end{proposition}

The following lemma explains why we can use definable $\kk$-partitions to compute the local density. It follows essentially from orthogonality between the value group and the residue rings. 
\begin{lemma}
\label{fibration}
Let $S$ and $A$ be definable, with $A\subseteq S[0,r,0]$. Then the following diagram commutes :

\[
\xymatrixcolsep{5pc}
\xymatrix{
\Cm_+(A[m,0,1]) \ar[d]^{\mu_{d,/A[0,0,1]}} \ar@/_3pc/[dd]_{\mu_{d,/S[0,0,1]}}&\\
\Cm_+(A[0,0,1]) \ar[d]^{\mu_{/S[0,0,1]}} \ar[r]^-{\MV} & \widetilde{\Cm_+}(A) \ar[d]^{\mu_{/S}} \\
\Cm_+(S[0,0,1]) \ar[r]_-{\MV} &\widetilde{\Cm_+}(S)}
\]

\end{lemma}
\begin{proof}
The left hand side commutes by definition of motivic integration. For the right hand side, one use the canonical isomorphism \[\Cm_+(A[0,0,1])\simeq \Qm_+(A)\otimes_{\Pp_+^0(S)}\Pp_+(S[0,0,1])\]
and the definition of $\MV$. 
\end{proof}

\subsection{Cones}

Kurdyka and Raby show in \cite{kurdyka_densite_1989} that the real local density of real subanalytic sets can be computed on the positive tangent cone. In our case, there is no canonical multiplicative subgroup of $K^\times$ anymore. One needs to consider a whole class of subgroups of $K^\times$, similarly to the $p$-adic case studied in \cite{cluckers_local_2012}.

\begin{definition}
Let $\Lambda$ a definable multiplicative subgroup of $K^\times$. A $\Lambda$-cone with origin $x$ is a definable $C\subseteq K^m$, with $x\in K^m$, such that for all $\lambda\in \Lambda$, and $y\in C$, 
\[ \lambda(y-x)\in C.\]

We say that $X\subseteq K^m$ is a local $\Lambda$-cone with origin $x$ if for some $\alpha\in \Zz$, $X=C\cap B(x,\alpha)$, with $C$ a $\Lambda$-cone with origin $x$. 

If $X\subseteq K^n$ and  $x\in K^n$, we call the definable set 
\[
\set{\lambda(y-x)\mid y\in X, \lambda \in \Lambda}
\]
the $\Lambda$-cone with origin $x$ generated by $X$. It is the smallest $\Lambda$-cone with origin $x$ that contains $X$.
\end{definition}

It will be sufficient to consider $\Lambda$-cones for $\Lambda \in \mathcal{D}$, where
\[
\mathcal{D}=\set{\Lambda_{n,m}\mid n,m\in \Nn^*},
\]
\[
\Lambda_{n,m}=\set{\lambda\in K^\times\mid \ac_m(\lambda)=1, \ord(\lambda)=0 \mod n}.
\]

\begin{remark}
\label{rem-cone}
\begin{enumerate}
\item If $X$ and $X'$ are respectively a $\Lambda_{n,m}$-cone and a $\Lambda_{n',m'}$-cone with origin $x$, then $X\cup Y$ is a $\Lambda_{\lcm(n,n'),\max(m,m')}$-cone with origin $x$. 

\item If $C$ is a $\Lambda_{n,m}$-cone with origin $x$, then 
\[
\Theta_d(C,x)=\frac{1}{n}\sum_{i=0}^{n-1}\frac{\mu_d(C\cap S(x,i))}{\eL^{-di}(1-\eL^{-d})}.
\]
\end{enumerate}
\end{remark}

\begin{definition}
Let $X\subseteq K^n$ be definable and $x\in K^n$. Let  
\[\pi_x : K^n\setminus \set{x}\to \Ppr^{n-1}(K)\] 
denote the function that send a point $z\neq x$ to the line between $z$ and $x$, and $\pi_x^X$ its restriction to $X\setminus\set{0}$. 
\end{definition}

We now show the local conic property of definable in ambient dimension 1. This is false in higher dimension, as is shown by the parabola example. 

\begin{lemma}
For every definable $X\subseteq K$, there is a $\Lambda\in \Dd$ and a definable function $\gamma : X \to \Zz$ such that $X\cap B(y,\gamma(y))$ is a local $\Lambda$-cone with origin $y$ for all $y\in Y$. 
\end{lemma}

\begin{proof}
By Remark \ref{rem-cone}, one can partition $X$ into finitely many definable parts. By cell decomposition, one can then assume that 
\[X=\set{x\in K \mid \ord(x-c)\in A, \ac_m(x-c)\in B},\]
 for a $c\in K$ and $A\subseteq \Zz$, $B\subseteq k_K$ definable. 
For $y\neq c$, we set $\gamma(y):=\ord(y-c)$, and $X\cap B(y,\gamma(y))$ is a local $\Lambda_{1,1}$-cone. For $y=c$, by Presburger quantifier elimination, up to cutting $A$ into a finite number of parts, there is a $b\in \Zz$ definable, $e\in \Nn^*$, $i\in\set{0,...,e-1}$ such that either
\[A\cap \set{\alpha\mid \alpha \geq b}=\set{\alpha\mid \alpha \geq b, \alpha=i \mod e},\]
 either this set is empty. In both cases, we set $\gamma(y)=b$, and $X\cap B(y,\gamma(y))$ is a local $\Lambda_{e,m}$-cone. 
\end{proof}

By working with parameters and using logical compactness, we get the following corollary. 

\begin{corollary}
\label{cor-lambda}
Let $X$ be a definable of $K^n$ and $x\in K^n$. Then there is a $\Lambda\in \Dd$ and a definable function $\alpha_x : \mathbb{P}^{n-1}(K) \to \Zz$ (\emph{i.e.} definable in affine charts of $\mathbb{P}^{n-1}(K)$) such that
\[
(\pi_x^X)^{-1}(\ell)\cap B(x,\alpha_x(\ell))
\]
is a local $\Lambda$-cone with origin $x$ for all $\ell\in \Pp^{n-1}(K)$. Moreover, one can choose $\Lambda$ independent from $x$, such that $(x,\ell)\mapsto \alpha_x(\ell)$ is definable. 
\end{corollary}
We call a $\Lambda\in \Dd$ satisfying the first condition adapted to $(X,x)$. If $\Lambda$ is adapted to $(X,x)$ for all $x$, then we say that $\Lambda$ is adapted to $X$.

\subsection{Tangent cone}
\label{subsection_tangent_cone}

\begin{definition}
\label{def-tan-con}
Let $\Lambda\in\Dd$, $X\subseteq K^n$ be definable and $x\in K^n$. The tangent $\Lambda$-cone of $X$ at $x$ is
\[
C_x^\Lambda(X)=\set{u\in K^n\mid \forall i\in \Zz, \exists y\in X, \exists \lambda\in \Lambda, \ord(y-x)\geq i, \ord(\lambda(y-x)-u)\geq i}.
\]
By definition, $C_x^\Lambda(X)$ is definable. 
\end{definition}
For X,Y definable in $K^n$, and $\Lambda'\subseteq \Lambda\in \Dd$, we have the following properties. 
\[
C_x^\Lambda(X\cup Y)=C_x^\Lambda(X)\cup C_x^\Lambda(Y),
\]
\[
C_x^\Lambda(\overline{X})=C_x^\Lambda(X),
\]
\[
C_x^{\Lambda'}(X)\subseteq C_x^{\Lambda}(X).
\]

\begin{remark}
\label{remark-tangent-cone}The inclusion above can be strict, but we have the following result. If $X$ is a local $\Lambda$-cone with origin $x$, then for all $\Lambda'\subseteq \Lambda$, $C_x^{\Lambda'}(X)= C_x^\Lambda(X)$. Indeed, as $X$ is a local $\Lambda$-cone, there is a $\Lambda$-cone $C$ with origin $x$ and an $n\in \Zz$ such that $X=C\cap B(x,n)$. We have $C_x^\Lambda(X)=C_x^\Lambda(C)$ and $C_x^{\Lambda'}(X)=C_x^{\Lambda'}(C)$. Let show the inclusion $C_x^{\Lambda}(C)\subseteq C_x^{\Lambda'}(C)$. Let $u\in C_x^{\Lambda}(C)$, there are sequences $(y_i)$ and $(\lambda_i)$   of elements of respectively $C$ and $\Lambda$ such that $y_i \to x$ and $\lambda_i(y_i-x) \to u$. As $\Lambda',\Lambda\in \Dd$, for all $i$ there are $\lambda_i'\in\Lambda'$, $\mu_i\in \Lambda$ such that $\abs{\mu_i}=1$ and $\lambda_i=\lambda_i'\cdot \mu_i$. Set $y_i'=\mu_i(y_i-x)+x$. Because $C$ is a $\Lambda$-cone $y_i'\in C$, and we have $y_i'\to x$, $\lambda_i'(y_i-x)\to u$ hence $u\in C_x^{\Lambda'}(C)$. 
\end{remark}

The following lemma justifies the use of $\kk$-partitions to compute the tangent cone, it follows from orthogonality between value group and residue rings. 
\begin{lemma}
\label{tang-con-kpart}
Let $X$ be a definable of $K^m$, and $g : X \to A$ a $\kk$-partition of $X$. 
Then for any $x\in K^m$ and $\Lambda \in \Dd$, 
\[
C_x^\Lambda(X)=\bigcup_{\xi\in A} C_x^\Lambda(X_\xi), 
\]
where $X_\xi=g^{-1}(\xi)$. 
\end{lemma}

\begin{proof}
Fix $u\in C_x^\Lambda(X)$. Let $D$ be the $u$-definable
\[
\set{(n,\xi)\in \Zz\times A\mid \exists y \in X, \exists \lambda \in \Lambda, \ord(y-x)\geq n, \ord(\lambda(y-x)-u)\geq n, g(y)=\xi}.
\]
By definition of $C_x^\Lambda(X)$, $D$ projects surjectively on its first coordinate. Because the theory splits, $D$ is a union of products of definable sets of $\Zz$ and $A$, hence there is a $\xi\in A$ such that for all $n\in \Zz$, $(n,\xi)\in D$. This means that $u\in C_x^\Lambda(X_\xi)$. 
\end{proof}

\begin{proposition}
\label{cone-spawned}
Let $X\subseteq K^m$ be definable, $x\in K^m$, and $\Lambda'\subseteq \Lambda \in \Dd$. Then $C_x^\Lambda(X)$ is the $\Lambda$-cone with origin 0 generated by $C_x^{\Lambda'}(X)$. 
\end{proposition}

\begin{proof}
Suppose $\Lambda'=\Lambda_{n,e}$. Consider the $\kk$-partition of $X$ defined by 
\[
g : X \to R_{e,K}^m, y\mapsto \ac_e(y-x).
\]
Let $u\in C_x^\Lambda(X)$ be nonzero. By Lemma \ref{tang-con-kpart} above, we can assume $u\in C_x^\Lambda(X_\xi)$, for some $\kk$-part $X_\xi$. Then there are some sequences $(y_k)_{k\geq 0}\in X_\xi$ and $(\lambda_k)_{k\geq 0}\in \Lambda$ such that $(y_k)$ converges to $x$ and $(\lambda_k(y_k-x))$ converges to $u$. By definition of the $\kk$-partition, $\ac_e(y_k-x)$ is constant. Because  $(\lambda_k(y_k-x))$ converges to $u$ and $u\neq 0$, $\ac_e(\lambda_k(y_k-x))$ is also constant for $k$ big enough. Hence $\ac_e(\lambda_k)$ is constant for $k$ big enough, so we can write (for $k$ big enough) $\lambda_k=\mu \lambda'_k$ where $\mu\in \Lambda$ does not depend on $k$ and $\ac_e(\lambda'_k)=1$. Up to extracting a subsequence, we can also assume that $\ord(\lambda'_k) \mod n$ is constant, hence up to multiplying $\mu$ and $\lambda'_k$ by some constant, we have $\lambda_k=\mu \lambda'_k$ where $\mu\in \Lambda$ and $\lambda'_k\in \Lambda'$. Then $\lambda'_k(y_k-x)$ converges to some $u'\in C_x^{\Lambda'}(X)$, and $u=\mu u'$. 
\end{proof}

\begin{corollary}
\label{dim-tang-cone}
The dimension of $C_x^\Lambda(X)$ is independent of $\Lambda\in \Dd$. 
\end{corollary}

\begin{proof}
As every $\Lambda\in \Dd$ is included in $\Lambda_{1,1}$, we can assume we are in the situation of Proposition \ref{cone-spawned}, with $\Lambda'\subseteq \Lambda\in \Dd$, and we have to prove that 
\[
\dim(C_x^{\Lambda}(X))=\dim(C_x^{\Lambda'}(X)).
\]
By Proposition \ref{cone-spawned}, $C_x^{\Lambda}(X)$ is the $\Lambda$-cone generated by $C_x^{\Lambda'}(X)$. Suppose $\Lambda=\Lambda_{n,e}$ and $\Lambda'=\Lambda_{n',e'}$. Then $n$ divides $n'$ and let 
\[
A=\set{\xi\in R_{e',K}\mid p_{e',e}(\xi)=1}\times \set{0,n,2n,...,n'-n}\subseteq R_{e',K}\times \Zz
\]
and 
\begin{eqnarray*}
B&=&\{(\xi,\alpha,z)\in A\times K^m \mid \exists y\in C_x^{\Lambda'}(X), \exists\mu\in \Lambda,\\ 
&&\hspace{1cm} \ac_e(\mu)=\xi, \ord(\mu)=\alpha, z=\mu y\}.
\end{eqnarray*}
Because $C_x^{\Lambda'}(X)$ is a $\Lambda'$-cone, the fiber $B_{(\xi,\alpha)}$ of $B$ above $(\xi,\alpha)\in A$ is in bijection with $C_x^{\Lambda'}(X)$. It follows that $D$ is of dimension $\dim(C_x^{\Lambda'}(X))$, because $A$ lives in the auxiliary sorts, hence is of dimension 0. Moreover, as $C_x^{\Lambda}(X)$ is the $\Lambda$-cone generated by $C_x^{\Lambda'}(X)$, the projection of $D$ on $K^m$ has image $C_x^{\Lambda}(X)$, hence 
\[
\dim(C_x^{\Lambda'}(X))\geq\dim(C_x^{\Lambda}(X)).
\]
The other inequality being clear by the inclusion $C_x^{\Lambda'}(X)\subseteq C_x^{\Lambda}(X)$, the corollary is proven.
\end{proof}

\begin{definition}
\label{def-deform}
Let $X$ a definable of $K^n$, $x\in K^n$ and $\Lambda\in \Dd$.  The deformation of $X$ is
\[
\Dd(X,x,\Lambda):=\set{(y,\lambda)\in K^n\times \Lambda \mid x+\lambda y\in X}
\]
Let $\overline{\Dd(X,x,\Lambda)}$ be its closure in $K^{n+1}$, we have $\overline{\Dd(X,x,\Lambda)}\cap K^n\times \set{0}=C_x^\Lambda(X)\times\set{0}$. We identify the latter with $C_x^\Lambda(X)$.
\end{definition}

\begin{proposition}
For every definable set $X$ of $K^n$, $x\in K^n$ and $\Lambda\in \Dd$, we have 
\[\dim(\Dd(X,x,\Lambda))=\dim(X)+1\quad \mathrm{\it{and}}\quad\dim(C_x^\Lambda(X))\leq \dim(X).\]
\end{proposition}

\begin{proof}
The projection $(y,\lambda)\in \Dd(X,x,\Lambda)\mapsto x+\lambda y\in X$ has its fibers isomorphic to $\Lambda$, hence of dimension 1. By dimension theory, $\dim(\Dd(X,x,\Lambda))=\dim(X)+1$. As $C_x^\Lambda(X)\subseteq \overline{\Dd(X,x,\Lambda)}\setminus \Dd(X,x,\Lambda)$, it is of dimension less than $\dim(\Dd(X,x,\Lambda))$. 
\end{proof}

We need to define a notion of multiplicity on the tangent cone. As shown by the following example, the multiplicity cannot be an integer. 
\begin{example}
Consider $K=\Cc((t))$, and the cusp $X=\set{y^2=x^3}$. The density of $X$ at 0 is 1/2. The tangent $\Lambda_2$-cone is

\begin{eqnarray*}
C_0^{\Lambda_2}(X)&=&\set{(0,y)\in K^2\mid \ord(y)=0\mod 2, \exists\eta\in \Cc^\times, \ac(y)^3=\eta^2}\\
&=&\set{(0,y)\in K^2\mid \ord(y)=0\mod 2},
\end{eqnarray*}

because $\Cc$ is algebraically closed. It seems that there are two fibers over each point $(0,y)$ of the tangent cone, defined by the two roots of $\ac(y)^3=\eta^2$, hence it seems we should attribute a multiplicity of 2 at every point of the tangent cone. But this is too naive, as it leads to a multiplicity of 1. A way to proceed is to view the multiplicity as a motive, the density $m_y$ of the tangent cone at $(0,y)$ is 
\[
[\set{\eta\in \Cc\mid \eta^2=\ac(y)^3}]\in \Cm_+(\set{y})
\]
and the tangent cone with multiplicities is
\[
CM_0^{\Lambda_2}(X)=[\set{(0,y,\eta)\in K^2\times\Cc\mid \ord(y)=0 \mod 2, \eta^2=\ac(y)^3}]\in \Cm_+(K^2),
\]
this leads to the correct value
\[
\Theta_1(CM_0^{\Lambda_2}(X))(0)=\frac{1}{2}\frac{[\{(\xi,\eta)\in k^2\mid \xi^2=\eta^3\}]-1}{\eL-1}=\frac{1}{2}.
\]
\end{example}
This motivates the following definition.

\begin{definition}
\label{def-cone-mult}
Let $X\subseteq K^n$ be definable of dimension $d$. Let $g : X \to A$ a definable $\kk$-partition of $X$ as in Proposition \ref{partition-1lipsch}. The definable $X$ is then partitioned into graphs $X_\xi$ of 1-Lipschitz functions defined on open sets of $K^d$, up to a set of dimension $<d$. 
For any $\Lambda\in \Dd$, we define the tangent $\Lambda$-cone with multiplicities as
\[
CM_x^\Lambda(X):= [\set{(y,\xi)\in K^n\times A \mid y\in C_x^\Lambda(g^{-1}(\xi))}]_d\in \mathrm{C}_+^d(K^n).
\]
\end{definition}

\begin{remark}
This depends a priori of the chosen partition. However, if we have two such $\kk$-partitions, by taking a refinement, we see that the class of 
\[m_y=\set{\xi\mid y\in C_x^\Lambda(g^{-1}(\xi))}\]
 does not depend on the chosen partition, up to a set of dimension $<d$. Hence $CM_x^\Lambda(X)$ is well defined (recall that $\mathcal{F}$unctions in $\mathrm{C}_+^d(K^n)$ are defined up to a set of dimension $<d$).
\end{remark}

We can now state our main theorem, that will be proved in Section \ref{section_proof_main_theorem}. 
\begin{theorem}
\label{thm-princ}
Let $X\subseteq h[n,0,0]$ a definable subassignment of dimension $d$. Then there is a $\Lambda\in \Dd$ such that for all $x\in K^d$
\[
\iota(\Theta_d(X,x))=\iota(\Theta_d(CM_x^\Lambda(X),0))\in \widetilde{\Cm}(\set{x}).
\]
\end{theorem}
Because the map from the Grothendieck semi-group to the group is not expected to be injective, the canonical induced map \[
\iota : \Cm_+(X)\to \Cm(X)
\]
is not expected to be injective either. As we will prove Theorem \ref{thm-princ} by a double inequality, this gives the equality only at the level of the Grothendieck group. This should not be seen as a problem, because every known realization of the Grothendieck semi-group factorizes through the Grothendieck group.

\section{Regular stratifications}
\label{section_regular_stratifications}

Regular stratifications satisfying the Whitney condition $(a_f)$ have been introduced by Thom \cite{thom_ensembles_1969} and studied by Hironaka \cite{hironaka_stratification_1977} in the complex case, Kurdyka-Raby \cite{kurdyka_densite_1989} in the real case. Afterward, the stronger Verdier condition $(w_f)$ has been introduced by Henry-Merle-Sabbah \cite{henry_sur_1984}, and studied by Bekka \cite{bekka_regular_1993} and Kurdyka-Parusinski \cite{kurdyka_$wsb_1994} in the real case, Loi \cite{loi_verdier_1998} in the $o$-minimal case and Cluckers-Comte-Loeser \cite{cluckers_local_2012} in the $p$-adic setting. This section is devoted to defining regular stratification satisfying analogous of conditions $(a_f)$ and $(w_f)$ in our context of tame Henselian fields. 

Note that Halupczok \cite{halupczok_non-archimedean_????} has already defined a notion of regular stratification in Henselian fields of equicharacteristic zero, the so-called $t$-stratifications. Instead of starting from the classical definition of $(a_f)$-stratification as we do, he starts by the property of local trivialization. This leads to two different notions of regular stratification in this context, and it is unknown whether a common generalization can be found, see \cite{halupczok_non-archimedean_????}, open question 9.1.

Because we work in an abstract (mixed) tame theory, we do not have a notion of analytic manifold, hence we only require in the following definition the strata to be differentiable. However, if one works in a concrete tame theory where a notion of analytic manifold is defined and an analytic cell decomposition theorem is available, then one can require the strata to be also analytic. 

\begin{definition}
A definable stratification of $Y\subseteq K^n$ is a definable $\kk$-partition of $Y$ such that each $\kk$-part is in definably diffeomorphic to some open subset of $K^d$ for some $d$ (depending on the chosen $\kk$-part). The $\kk$-parts are called strata and required to satisfy the so-called frontier condition :
\begin{center}
for each stratum $X$, $\overline{X}\setminus X$ is a union of strata.
\end{center}
If $f : Y \to K$ is a definable function, an $f$-stratification of $Y$ is a stratification of $Y$ such that for every stratum $X$, $f_{\vert X}$ is differentiable of constant rank. 
\end{definition}

Denote $T_{x,f}(X)=T_x(f_{\vert X}^{-1}(f(x)))$ the tangent space in $x$ of the fiber of $f_{\vert X}$ above $f(x)$. 

If $V,V'$ are two $K$-vector spaces of finite dimension with $\dim(V)\leq \dim(V')$, set
\[
\delta(V,V')=\sup_{v\in V,\abs{v}=1}\{\inf_{v'\in V',\abs{v'}=1} \abs{v-v'}\}.
\]
It is the Hausdorff distance between the units balls of $V$ and $V'$. 

It satisfies $\delta(V,V')=0$ if and only if $V\subseteq V'$, and for $V''\subseteq V'$ with $\dim(V'')\geq\dim(V)$, $\delta(V,V'')\geq\delta(V,V')$.

A couple of strata $(X,X')$ is said to satisfy the condition $(a_f)$ at $x_0 \in X\subseteq \overline{X'}$ if for every sequence $(x_n)$ of points of $X'$ converging to $x_0$, 
\[
\delta(T_{x_0,f} X, T_{x_n,f} X')\to 0 \tag{$a_f$}.
\]

A couple of strata (X,X') is said to satisfy the condition $(w_f)$ at $x_0 \in X\subseteq \overline{X'}$ if there is a constant $C$ such that for all neighborhoods $W$ of $x_0$, all $x\in W\cap X$ and $x'\in W\cap X'$, one has 
 \[
\delta(T_{x,f} X, T_{x',f} X')\leq C\abs{x-y} \tag{$w_f$}.
\]
A stratification $Y$ is said to be $(w_f)$-regular if for all strata $(X,X')$ and $x_0\in X\subseteq \overline{X'}$, the condition $(w_f)$ is satisfied at $x_0$.

If $V$ and $V'$ are definable, then the condition $\delta(V,V')\leq \epsilon$ is definable by the formula
 \[\forall v\in V, \abs{v}=1,\exists v'\in V', \abs{v'}=1,\abs{v-v'}\leq \epsilon.\]
 It follows that if $X,X'$ are definably diffeomorphic to some open subset of $K^d$ (for different $d$), then
 \[
 W_f(X,X')=\set{x\in X\mid (X,X') \mrm{satisfies} (w_f) \mrm{at} x}
\]
is definable. 

\begin{theorem}
\label{strat-reg}
Let $Y\subseteq K^n$ a definable, and $f : Y \to K$ a continuous definable function. Then there exists a $(w_f)$-regular definable stratification of $Y$. 
\end{theorem}

The proof of this theorem will follow closely the approach of Loi \cite{loi_verdier_1998} and uses as a key ingredient the curve selection lemma. We start by some lemmas. 

\begin{lemma}
\label{bound-diff}
Let $M : \Omega \to K$ be a differentiable and definable function on $\Omega \subseteq K^m\times K$.  Suppose $\overline{\Omega} \cap K^m\times \set{0}$ has non-empty interior $U\subseteq K^m$, and $M$ is bounded on $\Omega$. Then there exist a non-empty open $V\subseteq U$, $\epsilon>0$, and $d\in K^\times$ such that for all $(x,t)\in \Omega$, $x\in V$, and $\abs{t}\leq \epsilon$, one has
\[
\abs{D_xM(x,t)}\leq \abs{d}.
\]
\end{lemma}
Here and later, we use the notation $D_xM(x,t)=(\partial M(x,t)/\partial x_1, ... , \partial M(x,t)/\partial x_m)$.

\begin{proof}
By cell decomposition, we can partition $\Omega$ into $\kk$-definable parts such that on each cell $A$ such that $\overline{A}\cap K^m\times \set{0}\neq \emptyset$, one has
\[
\abs{D_{x_i}M(x,t)}=\abs{D_{x}M(x,t)}=\abs{h(x)}\cdot \abs{t}^a,
\]
for some $a\in \Qq$, $h$ a definable function, and $i\in \set{1,...,m}$. If all the $a$ are nonnegative, then we are done because the $\abs{h(x)}$ and the function bounding $\abs{t}$ from below are constant on small enough open subsets $V\subseteq U$. 

Assume that one $a$ is negative, say on a cell where $\abs{D_{x}M(x,t)}=\abs{\partial M(x,t)/\partial x_1}$. By the Jacobian property applied to $x_1\mapsto M(x_1,y)$, there is a  $(\kk,y)$-definable of $K$ such that the Jacobian property holds on each ball in it for the function $x_1\mapsto M(x_1,y)$. Up to another cell decomposition, one can assume that there are some $\kk$-definable balls $B_1\subseteq K$, $B\subseteq K^{m-1}$ and $C\subseteq K$ a $\kk$-definable with $0\in \overline{C}$ such that for all $y\in B\times C$, the Jacobian property holds on $B_1$ for $x_1\mapsto M(x_1,y)$. We then have for all $y=(x_2,...,x_n,t)\in B\times C$ and $x_1,x_1'\in B_1$,
\[
\abs{M(x_1',y)-M(x_1,y)}=\abs{x_1'-x_1}\cdot \abs{\partial M(x,t)/\partial x_1}=\abs{x_1'-x_1}\cdot \abs{h(x)}\cdot\abs{t}^a
\]
As $M$ is bounded, the left term is bounded but the right one goes to $+\infty$ when $t$ goes to $0$, this is a contradiction. 
\end{proof}

\begin{lemma}
\label{trivialization}
Let $V\subseteq K^m$, $S\subseteq K^n$ be definable with $V$ nonempty open in $K^m$, and $V\times\set{0}^{m-n}\subseteq \overline{S}\setminus S$. Let $\pi : K^n \to K^m$ be the natural projection, and for every $u\in V$, $S_u=\set{x\in S \mid \pi(x)=u}$. Then there is a nonempty open $U\subseteq V$ such that for all $u\in U$, $(u,0)\in \overline{S_u}\setminus S_u$.
\end{lemma}

\begin{proof}
Let \[A=\set{(x,y,\delta)\in S\times V\times \Gamma\mid \pi(x)=y, \abs{x-y}\leq \delta}.\]
Let $\pi_1$ the projection defined by $\pi_1(x,y,\delta)=y$. The definable $\pi_1(A)$ is dense in $V$. Indeed, let $y_0\in V$, and $\epsilon>0$. By openness of $V$, there is an $\epsilon'$ with $0<\epsilon'<\epsilon$ and $B(y_0,\epsilon')\cap K^m\times \set{0}^{n-m}\subseteq V\times \set{0}^{n-m}$. Since $y_0\in \overline{S}\setminus S$, there is a $x\in B(y_0,\epsilon')\cap S$. Then $y:=\pi(x)\in \pi_1(A)$ and $\abs{y-y_0}\leq \epsilon$. 

Let $\pi_2$ be the projection defined by $\pi_2(x,y,\delta)=(y,\delta)$ and set
\[
\pi_2(A)_y := \set{\delta\mid (y,\delta)\in \pi_2(A)},\; \epsilon(y)=\inf \pi_2(A)_y, \; y \in \pi_1(A)
\]
Then $\epsilon$ is definable, and it suffices to show that $\dim(\set{ y\in \pi_1(A)\mid \epsilon(y)>0})<m$. Suppose the contrary. Then by cellular decomposition there is a $W\subseteq \pi_1(A)$ open on which $\epsilon$ is greater than a $c>0$. 
But then for all $y_0\in W$, and $\delta$ with $0<\delta<c$ and $B(y_0,\delta)\cap K^m\times \set{0}^{n-m}\subseteq W\times \set{0}^{n-m}$. Then $\abs{x-y}>c$ for all $y\in B(y_0,\delta)\cap K^m$, $x\in S$ with $\pi(x)=y$. Contradiction with the argument above. 
\end{proof}

\begin{lemma}
\label{lemma-strat}
Let $f : Y\subseteq K^n \to K$ a definable continuous function. Assume that for all couples $(X,X')$ of definable differentiable sets, with $X\subseteq \overline{X'}\setminus X'$, we have $\dim(X\setminus W_f(X,X'))\leq \dim(X)$. Then $Y$ admits a $(w_f)$-regular $f$-stratification. 
\end{lemma}

\begin{proof}
There exists $S$ an $f$-stratification of $Y$ by Proposition \ref{partition-1lipsch} applied inductively. By decreasing induction on $d\in\set{0,...,n}$, one builds some successive refinements $S_d$ of $S$ satisfying :
\[
\mrm{for \; all} X,X'\in S_d, X\subset \overline{X'}, \dim(X)\geq d, \mrm{one \; has} W_f(X,X')=X,
\] 
using the lemma hypothesis. Then $S_0$ is ($w_f$)-regular. 
\end{proof}

\begin{proof}[Proof of Theorem \ref{strat-reg}]
By Lemma \ref{lemma-strat}, it suffices to show that for a pair $(X,X')$ of definable manifolds of $Y$, with $X\subseteq \overline{X'}\setminus X'$, one has $\dim(X\setminus W_f(X,X'))<\dim(X)$. Suppose by contradiction that for such a pair, $\dim(X\setminus W_f(X,X'))=\dim(X)$. Up to reducing $X$ and $X'$, assume $W_f(X,X')$ is empty, that $f$ is of constant rank on respectively $X$ and $X'$. We can also replace $X$ by $f_{\vert X}^{-1}(a)$, for a $a\in K$ such that this set is nonempty, therefore assume $f$ is constant on $X$.  By invariance of the condition $(w_f)$ under diffeomorphisms, one can assume $X=V\times \set{0}^{n-m}$, for $V\subseteq K^m$ open definable, therefore for $x\in X$, $T_x X=K^m\times\set{0}^{n-m}$.  

 There are two cases to consider. 

$\bullet$ {\bf Case 1 :}  $f$ is of rank $0$ on $X'$. In this case, one can assume $f$ constant on $X'$. As $W_f(X,X')=\emptyset$, the definable 
\[
D=\set{(y,t)\in X'\times K^\times\mid  \delta(K^m\times \{0\}^{n-m}, T_{y} X')\geq \abs{y_1} \abs{t}^{-1}}
\]
satisfies $V\times\set{0}^{n-m+1}\subseteq \overline{D}\setminus D$ (where $y=(y_0,y_1)$, with $y_1$ of length $n-m$). By Lemma \ref{trivialization}, there is a nonempty open $U\subseteq V$ such that for all $u\in U$, $(u,0)\in \overline{D_u}\setminus D_u$. By the curve selection lemma $\ref{CSL}$, there is a definable (with parameters) function $\overline{\rho} : U\times B_e \to \overline{D}$ with $\overline{\rho}(u,s)\in \overline{D}_u$, $\overline{\rho}(u,0)=(u,0)$, $\overline{\rho}(u,s)\in D$ if $s\neq 0$, $\overline{\rho}(u,\cdot)$ $C$-Lipschitz. By the Jacobian property, up to restricting ourself to an open subset of $U$ one can moreover assume that $\overline{\rho}(\cdot,s)$ is differentiable. Also denote $\overline{\rho}(u,s)=(u,\rho(u,s),r(u,s))$. We have for all $u\in U$ and $s\in B_e^*$, 
\[
\delta(K^m\times \{0\}^{n-m}, T_{u,\rho(u,s)} X')\geq \abs{\rho(u,s)} \abs{r(u,s)}^{-1}.
\]
By definition of Hausdorff distance, 
\[
\delta(K^m\times \{0\}^{n-m}, T_{u,\rho(u,s)} X')\leq \abs{D_u\rho(u,s)}.
\]
By cell decomposition, up to replacing $(u,s)\mapsto \overline{\rho}(u,s)$ by $(u,s)\mapsto \overline{\rho}(u,s^k)$, and $U$ by a nonempty open subset, there are an $a\in K^\times$ and an integer $\ell>0$ such that for all $(u,s)\in U\times B$,
\[
\abs{\rho(u,s)}=\abs{a}\abs{s}^\ell. 
\]
Since the function $(u,s)\in U\times B^* \to \rho(u,s)/s^\ell$ is bounded, by Lemma \ref{bound-diff}, there is a $d \in K^\times$ such that up to replacing $U$ by a nonempty open subset, for all $(u,s)\in U\times B^*$, 
\[
\abs{D_u\rho(u,s)}\leq \abs{d}, \mrm{hence} \frac{\abs{D_u\rho(u,s)}}{\abs{\rho(u,s)}}\leq  \abs{ad}.
\]
On the other hand, 
\[
\frac{\abs{D_u\rho(u,s)}}{\abs{\rho(u,s)}}\geq \abs{r(u,s)}^{-1},
\]
and $r(u,s)$ goes to $0$ when $s$ goes to $0$, contradiction.

$\bullet$ {\bf Case 2 :} $f$ is of rank 1 on $X'$. We can  assume $f(X)=0$ and $f(X')\neq 0$. As $W_f(X,X')=\emptyset$, the definable
\[
D=\set{(y,t)\in X'\times K^\times\mid  \delta(K^m\times \{0\}^{n-m}, T_{y,f} X')\geq \abs{y_1} \abs{t}^{-1}}
\]
satisfies $V\times \set{0}^{n-m+1}\subseteq \overline{D_u}\backslash D_u$ (where $y=(y_0,y_1)$, with $y_1$ of length $n-m$). By Lemma \ref{trivialization}, there is a nonempty open $U\subseteq V$ such that for all $u\in U$, $(u,0)\in \overline{D_u}\backslash D_u$. By Corollary \ref{cor-CSL}, up to restricting $U$ to a $\kk$-definable subpart, there is a definable with parameters function $\overline{\rho} : U \times B_e\to \overline{D}$ with $\overline{\rho}(u,s)\in \overline{D_u}$, $\overline{\rho}(u,0)=(u,0)$, $\overline{\rho}(u,s)\in D$ if $s\neq 0$ and $\overline{\rho}(u,\cdot)$ $C$-Lipschitz for all $u\in U$. Denote $\overline{\rho}(u,s)=(u,\rho(u,s),r(u,s))$. By Proposition \ref{inj-cons}, we can moreover assume (up to restricting again $U$ and taking $e$ bigger) that $f(u,{\rho}(u,\cdot))$ is injective (it cannot be constant because $f$ is of rank one on $X'$).   Consider $f_0(u,s)=f(u,\rho(u,s))$. By cell decomposition, there is an integer $e'$ and a definable with parameters continuous function $h : B_{e'} \to f_0(U\times B_e)$ with $h(0)=0$ and $h(B_{e'}^*)\subseteq K^\times$. For $(u,s')\in U\times B_{e'}$, there is an unique $s\in B_e$ such that $f_0(u,s)=h(s')$, denote $h_0(u,s')=s$. The function $h_0$ is definable (with parameters) and $h_0(u,s')$ goes to $0$ when $s'$ goes to 0 (for fixed $u$). Replace $\overline{\rho}$ by 
\[
\fonction{\overline{\rho'}}{U\times B_{e'}}{\overline{D}}{(u,s)}{\overline{\rho}(u,h_0(u,s))}
\]
The function $\overline{\rho'}$ satisfies $\overline{\rho'}(u,s)\in \overline{D_u}$, $\overline{\rho'}(u,0)=(u,0)$, $\overline{\rho'}(u,s)\in D$ if $s\neq 0$ and $\lim_{s\to 0}\overline{\rho'}(u,s)=0$ for all $u\in U$. Up to reducing $U$, assume $\overline{\rho'}(\cdot,s)$ is differentiable with respect to $u$ at all $u\in U$. Denote $\overline{\rho'}(u,s)=(u,\rho'(u,s),r'(u,s))$. For all $(u,s)\in U\times B_{e'}$, 
\[
\delta(K^m\times \{0\}^{n-m}, T_{(u,\rho(u,s)),f} X')\geq \abs{\rho'(u,s)} \abs{r'(u,s)}^{-1}
\]
and
\[
f(u,\rho'(u,s))=f_0(u,s)=h_0(s).
\]
The vector space spanned by vectors $(0,...,0,1,0,...,0, \partial \rho'(u,s)/\partial u_i)$ (with 1 in $i$-th position) is a subspace of $T_{(u,\rho'(u,s)),f} X'$, hence 
\[
\abs{D_u \rho'(u,s)}\geq \delta(K^m\times \{0\}^{n-m}, W)\geq \delta(K^m\times \{0\}^{n-m}, T_{(u,\rho(u,s)),f} X').
\]
We are now in a situation similar to case 1, up to replacing $\overline{\rho'}(u,s)$ by $\overline{\rho'}(u,s^k)$ and reducing $U$, we can assume 
\[
\abs{\rho'(u,s)}=\abs{a}\abs{s}^\ell,
\]
and by Lemma \ref{bound-diff}, 
\[
\abs{D_u\rho((u,s))}\leq \abs{d}.
\] 
Hence 
\[
\abs{r'(u,s)}^{-1}\leq \frac{\abs{D_u\rho(u,s)}}{\abs{\rho(u,s)}}\leq \abs{ad}, 
\]
contradiction.
\end{proof}

\section{Proof of Theorem \ref{thm-princ}}
\label{section_proof_main_theorem}
\subsection{Codimension 0 case}
The proof of Theorem \ref{thm-princ} goes along the same lines as the proof of the main theorem of \cite{cluckers_local_2012}. We start by proving it in the case $d=n$, using a monotone convergence lemma.

\begin{lemma}[Monotone convergence]
\label{mon-conv}
Let $\varphi \in \Cm_+(A\times \Nn)$ integrable. Then 
\[
\lim_{n\to +\infty}\sum_{k=0}^n \mu(i_n^*(\varphi))=\mu(\varphi).
\] 
In particular, if $(X_n)_{n\geq 0}$ is a bounded uniformly definable sequence of sets of $K^m$, increasing for  the inclusion (resp. decreasing), then
\[
\lim_{n\to \infty}\mu(\bigcup_{0\leq k\leq n}X_k)=\mu(\bigcup_{k=0}^{+\infty} X_k),
\]
\[
\mrm{(resp.)}\;\lim_{n\to \infty}\mu(\bigcap_{0\leq k\leq n}X_k)=\mu(\bigcap_{k=0}^{+\infty} X_k).
\]
\end{lemma}
\begin{proof}
The sequence $\sum_{k=0}^n \mu(i_n^*(\varphi))$ is bounded, definable, and increasing. By Presburger cell decomposition, as in Lemma \ref{lem-MVPP}, we have the existence of the limit. By definition of the motivic integral with respect to a value group variable, we have the result. 

We consider now the definable $X=\bigcup_{0\leq k} X_k$, where $(X_k)_k$ is an increasing sequence. Let $f : X \to \Nn$ the definable function that maps $x$ to the least $n$ such that $x\in X_n$. Let $\widetilde{X}_n=f^{-1}(n)=X_n\setminus X_{n-1}$. We see $X$ as a definable of $X\times\Nn$, and by the first point,
\[\lim_{n\to+\infty}\sum_{k=0}^n \mu(\widetilde{X}_k)=\mu(X)=\mu(\bigcup_{k=0}^{+\infty}X_k).\]
By linearity of the integral,
\[ \sum_{k=0}^n \mu(\widetilde{X}_k)=\mu(\bigcup_{0\leq k\leq n}\widetilde{X}_k)=\mu(\bigcup_{0\leq k\leq n}X_k),\]
hence we have the required equality ; the other case is similar. 
\end{proof}

\begin{proposition}
\label{full-dim}
Let $X$ be definable of $K^d$. Then for all $\Lambda\in \Dd$ small enough and $x\in K^d$,
\[
\iota(\Theta_d(X,x))=\iota(\Theta_d(C_x^\Lambda(X),0)).
\]
\end{proposition}
\begin{proof}
We have to show the two inequalities 
\[
\Theta_d(C_x^\Lambda(X),0)\geq \Theta_d(X,x)\]
and
\[
\Theta_d(C_x^\Lambda(X),0)\leq \Theta_d(X,x).
\] 

By Corollary \ref{cor-lambda}, there is a $\Lambda=\Lambda_{n,m}\in \Dd$, a definable function $\alpha : \Ppr^{d-1}(K)\to \Zz$ such that for all $\ell\in \Ppr^{d-1}(K), (\pi_x^X)^{-1}(\ell)\cap B(x,\alpha(\ell))$ is a local $\Lambda$-cone with origin $x$. This remains true for any $\Lambda'\in \Dd$, $\Lambda'\subseteq \Lambda$, hence the following reasoning holds also for any such $\Lambda'$. 

For all $n\geq 0$, let $C_n$ the $\Lambda$-cone with origin $x$ generated by $X\cap B(x,n)$ and $W=\bigcap_{n\geq 0}C_n$, $W$ is also a $\Lambda$-cone. By monotone convergence lemma \ref{mon-conv}, for any $i=0,...,n-1$
\[
\mu_d(W\cap S(x,i))=\lim_{n\to +\infty} \mu_d(C_n\cap S(x,i)),
\]
hence by Remark \ref{rem-cone} and linearity of the integral, 
\[
\Theta_d(W,x)=\lim_{n\to +\infty} \Theta_d(C_n,x).
\]
Moreover, the tangent cone satisfies by definition $x+C_x^\Lambda(X)=\cap_{n\geq 0}\overline{C_n}$. The $\overline{C_n}$ being also $\Lambda$-cones, for the same reason we have
\[
\Theta_d(C_x^\Lambda(X),0)=\lim_{n\to +\infty} \Theta_d(\overline{C_n},x).
\]
For all $n$ we have $\Theta_d(\overline{C_n},x)=\Theta_d(C_n,x)$, so
\[
\Theta_d(W,x)=\Theta_d(C_x^\Lambda(X),0). 
\]
Because $\Theta_d$ is increasing, for all $n\geq 0$, we have
\[
\Theta_d(X,x)=\Theta_d(X\cap B(x,n))\leq \Theta_d(C_n,x), 
\]
hence $\Theta_d(X,x)\leq \Theta_d(C_x^\Lambda(X),0)$.

For the other inequality, for all $n\geq 0$, set $W_n$ defined as the $w\in W$ such that $\alpha(\pi_x^X(w))\leq n$. By definition of $\alpha$, $\bigcup_{n\geq 0} W_n=W\setminus\set{0}$, and $W_n\cap B(x,n)$ is a local $\Lambda$-cone. As before, by the monotone convergence lemma \ref{mon-conv}
\[\Theta_d(W,x)=\lim_{n\to +\infty} \Theta_d(W_n,x).\]
Moreover, for all $w\in W_n\cap B(x,n)$, there is a $y\in X\cap B(x,n)$ and a $\lambda\in \Lambda$ such that $w=\lambda(y-x)+x$ and $\alpha(\pi_x^X(w))\leq n$, so $w\in X$. We showed that $W_n\cap B(x,n)\subseteq X\cap B(x,n)$, hence
\[
\Theta_d(W_n,x)\leq \Theta_d(X,x), 
\]
and by taking the limit, 
\[
\Theta_d(C_x^\Lambda(X),0)=\Theta_d(W,x)\leq \Theta_d(X,x).
\]
\end{proof}

\subsection{Deformation to the tangent cone.} The above proof does not extend to the case $X\subseteq K^n$ of dimension $d$, for $d<n$. The reason is that the $\Lambda$-cone generated by a definable of dimension $d$ is in general of dimension $n$. However, using Proposition \ref{partition-1lipsch}, we will restrict ourself to the case of a graph of a 1-Lipschitz function of domain $U$ included in $K^d$. By the case we just studied, it suffices to show that for some $\Lambda$, \[\Theta_d(C_u^\Lambda(U),0)=\Theta_d(C_x^\Lambda(X),0).\]
To do so, using deformation to the tangent cone we will build a 1-Lipschitz function $\psi$ of domain (included in) $C_u^\Lambda(U)$ such that the graph of $\psi$ is a dense subset of $C_x^\Lambda(X)$. The displayed equality follows from the existence of such a $\psi$. We need first some preparatory lemmas.

\begin{remark}
\label{rem-lambda-unique}
In the following lemmas, at several occasions we will have to use a $\kk$-partition of the studied definable. By applying Corollary \ref{cor-lambda} with parameters, we can assume that the same $\Lambda$ is adapted to all the $\kk$-parts. 
\end{remark}

Let $U$ an open of $K^d$. By Corollary \ref{cor-lambda}, there is a definable function 
\[\alpha : \Ppr^{d-1}(K) \to \Zz\cup\set{\infty}\]
such that for all $\ell\in\Ppr^{d-1}(K)$, $(\pi_u^U)^{-1}(\ell)\cap B(u,\alpha(\ell))$ is a local $\Lambda$-cone with origin $u$, with the convention that $\alpha(\ell)=\infty$ if and only if $(\pi_u^U)^{-1}(\ell)=\emptyset$. Because $\alpha$ is definable, it is continuous on a definable dense open subset $\Omega_0$ of $\Ppr^{d-1}(K)$. By dense, we mean that the dimension of its complement into $\Ppr^{d-1}(K)$ is strictly less than $d-1$. Let fix an $u\in K^d$ and let $\Omega_1$ the definable constituted by the $\ell\in\Omega_0$ such that for all neighborhood $V$ of $u$ in $K^d$, $(\pi_u^U)^{-1}(\ell)\cap V$ is nonempty. 

\begin{lemma}
\label{lemma-dens1}
Suppose $C_u^\Lambda(U)$ is of maximal dimension, \emph{i.e.} $d$. Then $\Omega_1$ contains a nonempty open definable of $\Ppr^{d-1}(K)$. 
\end{lemma}

\begin{proof}
Let $\Omega_1^c$ the complement of $\Omega_1$ in $\Ppr^{d-1}(K)$. We need to show that $\Omega_1^c$ is not dense. Assume it is, for the sake of contradiction. By definability and density, $\Omega_1^c\cap \Omega_0$ is also dense in $\Ppr^{d-1}(K)$. For every $\ell\in \Omega_1^c\cap \Omega_0$, by definition of the tangent cone, $(\pi_u^U)^{-1}(\ell)\cap (u+C_u^\Lambda(U))=\emptyset$. As $(\pi_u^U)^{-1}(\Omega_1^c\cap \Omega_0)$ is dense and definable in $K^d$, this implies that $C_u^\Lambda(U)$ is contained into a definable of dimension less than $d$, contradiction. 
\end{proof}

For all definable $O$ of $\Ppr^{d-1}(K)$, we consider the definable
\[
C_u^{\Lambda,O}(U):=(\pi_u)^{-1}(O)\cap C_u^\Lambda(U).
\]
\begin{lemma}
\label{lemma-dens2}
We still suppose $\dim(C_u^\Lambda(U))=d$. Let $O$ a dense open definable of $\Omega_1$. Then $C_u^{\Lambda,O}(U)$ is dense open in $C_u^\Lambda(U)$, hence of dimension $d$. 
\end{lemma}

\begin{proof}
As $\pi_u$ is continuous, $(\pi_u)^{-1}(O)$ is open, so $C_u^{\Lambda,O}(U)$ is open in $C_u^\Lambda(U)$. Let show that it is dense. Its complement is the disjoint union 
\[
(C_u^\Lambda(U)\cap (\pi_u)^{-1}(\Omega_1\setminus O))\cup (C_u^\Lambda(U)\cap (\pi_u)^{-1}(\Omega_1^c)).
\]
By density of $O$ in $\Omega_1$, $C_u^\Lambda(U)\cap (\pi_u)^{-1}(\Omega_1\setminus O)$ is of dimension less than $d$. For all $z\in C_u^\Lambda(U)$, $\pi_u(z)\in \overline{\Omega_1}$, hence
\[
C_u^\Lambda(U)\cap (\pi_u)^{-1}(\Omega_1^c)\subseteq (\pi_u)^{-1}(\overline{\Omega_1}\setminus\Omega_1),
\]
this shows that it is of dimension less that $d$. 
\end{proof}

\begin{lemma}
\label{lemma-equal}
For almost all $\ell\in \Omega_1$ (\emph{i.e.} in a dense subset of $\Omega_1$), 
\[
C_u^\Lambda(U)\cap \ell \cap B(u,\alpha(\ell))=\overline{(\pi_u^U)^{-1}(\ell)} \cap B(u,\alpha(\ell)).
\]
\end{lemma}

\begin{proof}
Suppose $u=0$ to simplify the notations. For any $\ell\in \Omega_1$, we have
\begin{eqnarray*}
\overline{(\pi_u^U)^{-1}(\ell)}\cap B(u,\alpha(\ell))&=&C_u^\Lambda((\pi_u^U)^{-1}(\ell))\cap B(0,\alpha(\ell))\\
&\subseteq&C_u^\Lambda(U)\cap \ell \cap B(u,\alpha(\ell)), 
\end{eqnarray*}
hence one inclusion is proven. For the other one, let $X$ be the definable set
\[
X:=\set{(x,t)\in K^d\setminus\set{0}\times K\mid u+tx\in U}, 
\]
parametrizing lines $\ell\cap U$. Let us apply a cell decomposition of $X$ relatively to the $t$ variable. We recall the notation $X_x$ for the fiber of $X$ above $x\in K^d\setminus\set{0}$. There are three cases. Either 0 is in the interior of $X_x$, or is in its (definable) boundary $\partial X_x$, or it is outside $\overline{X_x}$. We do not need to consider the last case, because $\ell\in \Omega_1$. For the first one, $(\pi_u^U)^{-1}(\ell)\cap B(u,\alpha(\ell))=B(u,\alpha(\ell))$, hence the inclusion is satisfied. For the case $0\in \partial X_x$, up to a set of $\ell$ of dimension less than $d-1$, we also have the inclusion, because the function (of the variable $x$) describing the centers of the cell is almost everywhere continuous. 
\end{proof}

\begin{lemma}
\label{lemma-omega}
We suppose $\dim(C_u^\Lambda(U))=d$. Then there is a dense open definable subset of $\Omega_1$ such that for all $z\in C_u^{\Lambda,\Omega}(U)$, with $\ell:=\pi_u(z)$, and all $\lambda\in \Lambda$ small enough (depending on $z$),
\[
u+\lambda z \in (\pi_u^U)^{-1}(\ell).
\]
\end{lemma}

\begin{proof}
For all $\ell\in \Omega_1$, $z\in \ell \cap C_u^\Lambda(U)$, and $\lambda\in \Lambda$, we have $u+\lambda z\in C_u^\Lambda(U)$.  From the inclusion $\subseteq$ of Lemma \ref{lemma-equal}, we have that for all $\lambda\in \Lambda$ with $\ord(\lambda z)\geq \alpha(\ell)$,
\[u+\lambda z \in \overline{(\pi_u^U)^{-1}(\ell)}.\]
As $\overline{(\pi_u^U)^{-1}(\ell)}\setminus (\pi_u^U)^{-1}(\ell)$ is finite, we have for every $\lambda\in\Lambda$ small enough $u+\lambda z \in (\pi_u^U)^{-1}(\ell)$. 
\end{proof}

\begin{corollary}
\label{cor-stab-cone}
Let $d>0$, and $U$ a definable nonempty open of $K^d$, and $\Lambda\in\Dd$ as in Corollary \ref{cor-lambda}. Then for all $u\in \overline{U}$ such that $C_u^\Lambda(U)$ is of maximal dimension $d$, for every $\Lambda'\in \Dd$ with $\Lambda'\subseteq \Lambda$, we have 
\[C_u^{\Lambda'}(U)=C_u^\Lambda(U).\]
\end{corollary}

\begin{proof}
Let us show the inclusion $C_u^{\Lambda,\Omega}(U)\subseteq C_u^{\Lambda',\Omega}(U)$, with $\Omega$ given by Lemma \ref{lemma-omega}. Let $z\in C_u^{\Lambda,\Omega}(U)$, and $\ell=\pi_u(z)$. We have $z\in C_u^\Lambda(\overline{U\cap \ell})=C_u^\Lambda(U\cap \ell)$, hence by Lemma \ref{lemma-omega}, $U\cap \ell \cap B(u,\alpha(\ell))$ is a local $\Lambda$-cone with origin $u$, hence by Remark \ref{remark-tangent-cone},
\[
z\in C_u^\Lambda(U\cap \ell)=C_u^{\Lambda'}(U\cap \ell)\subseteq C_u^{\Lambda'}(U).
\] 
The other inclusion being clear, we showed that $C_u^{\Lambda,\Omega}(U)= C_u^{\Lambda',\Omega}(U)$, hence that $\overline{C_u^{\Lambda,\Omega}(U)}= \overline{C_u^{\Lambda',\Omega}(U)}$. By Lemma \ref{lemma-dens2}, $\overline{C_u^{\Lambda,\Omega}(U)}=C_u^\Lambda(U)$. But Lemmas \ref{lemma-dens1} and \ref{lemma-dens2} are also true for $\Lambda'$, provided that $\dim(C_u^{\Lambda'}(U))=\dim(C_u^{\Lambda}(U))=d$ by Corollary \ref{dim-tang-cone}. Hence
\[
C_u^{\Lambda'}(U)=\overline{C_u^{\Lambda',\Omega}(U)}=\overline{C_u^{\Lambda,\Omega}(U)}=C_u^\Lambda(U).
\]
\end{proof}

\begin{proposition}
\label{prop-psi1}
Let $U$ an open definable of $K^d$, and $\varphi : U \to K^{m-d}$ definable, differentiable and locally 1-Lipschitz. Let $u\in \overline{U}$ be fixed, $\Lambda\in \Dd$ adapted to $(U,u)$. We suppose $\lim_{x\to u} \varphi(x)=v$, and set $w=(u,v)$. We suppose $\dim(C_u^\Lambda(U))=d$, and consider $\Omega$ and $C_u^{\Lambda,\Omega}(U)$ as given by Lemma \ref{lemma-omega}.

Then there exists a $\kk$-partition of $U$ and $\Lambda'\in \Dd$, such that for each $\kk$-part $U_\xi$ with $\dim(U_\xi)=d$, $\varphi$ is globally $C$-Lipschitz on $U_\xi$, $\Lambda'$ is adapted to $(U_\xi,u)$, there is a $\Omega_\xi$ given by Lemma \ref{lemma-omega}. For all $z\in C_u^{\Lambda',\Omega_\xi}(U_\xi)$, the limit 
\[
\underset{\lambda\in \Lambda'}{\lim_{\lambda \to 0}} \;\lambda^{-1}(\varphi(u+\lambda z)-v)
\]
exists, yielding a definable function $\psi : C_u^{\Lambda',\Omega_\xi}(U_\xi) \to K^{m-d}$.
\end{proposition}

\begin{proof}
By Lemma \ref{lemma-omega}, for all $z\in C_u^{\Lambda,\Omega}(U)$, for all $\lambda\in\Lambda$ small enough, $u+\lambda z\in U$. Define $U'\subseteq U$ as the set of all such $u+\lambda z$. By definition, $C_u^{\Lambda,\Omega}(U)\subseteq C_u^\Lambda(U')$, hence $\dim(C_u^\Lambda(U\setminus U'))<d$, we can neglect it and suppose $U'=U$. By Theorem \ref{thm-global-lipsch}, assume that $\varphi$ is globally $C$-lispchitz on $\kk$-parts of a $\kk$-partition of $U$. The function 
\[
f_z : \lambda \mapsto \frac{\varphi(u+\lambda z)-v}{\lambda}
\]
is bounded, indeed otherwise there would exist a $y\in K^{m-d}\backslash\set{0}$ such that $(0,y)\in C_u^\Lambda(\Gamma(\varphi))$, in contradiction with Proposition \ref{prop-tangent-incl} and the fact that $\varphi$ is $C$-Lipschitz. By Lemma \ref{lemma-adh-value}, $f_z$ has at least one limit value when $\lambda \to 0$. 
As $\dim(\overline{\Gamma(f_z)}\setminus \Gamma(f_z))=0$, this set is finite, so the $z$-definable $A_z\subseteq K^{m-d}$ of adherence values at 0 of $f_z$ is nonempty and finite. By Lemma \ref{lemma-finite-bij}, there is a $z$-definable function $g_z : K^{m-d} \to R_{s,K}^s$ sending bijectively $A_z$ to some $B_z$ and we can assume $g_z$ is constant on a neighborhood of each point of $A_z$. We can moreover assume that if $z'=\lambda' z$ for some $\lambda'\in\Lambda$, then $g_{z'}(a)=g_z(a/\lambda')$. The definable function 
\[
u+\lambda z \in U \mapsto g_z\left(\frac{\varphi(u+\lambda z)-v}{\lambda}\right)
\]
is then well defined, and yields a $\kk$-partition of $U$. We can then apply Corollary \ref{cor-lambda}, Remark \ref{rem-lambda-unique}, and Lemmas \ref{lemma-dens1} to \ref{lemma-omega} to get $\Lambda'\subseteq \Lambda$ and $\Omega_\xi$ as required. For any $z\in C_u^{\Lambda',\Omega_\xi}(U_\xi)$ and $\lambda\in \Lambda'$ small enough such that $u+\lambda z\in U_\xi$, we have $g_z(\frac{\varphi(u+\lambda z)-v}{\lambda})$ constant, so $\frac{\varphi(u+\lambda z)-v}{\lambda}$ has a limit when $\lambda\in \Lambda'$ goes to 0. 
\end{proof}

\begin{proposition}
\label{prop-psi2}
Let $U$ be an open definable of $K^d$, and $\varphi : U \to K^{m-d}$ definable, differentiable and locally 1-Lipschitz and globally C-Lipschitz. Let $u\in \overline{U}$ be fixed, $\Lambda\in \Dd$ adapted to $(U,u)$. We suppose $\lim_{x\to u} \varphi(x)=v$, and note $w:=(u,v)$. We suppose $\dim(C_u^\Lambda(U))=d$, and consider $\Omega$ and $C_u^{\Lambda,\Omega}(U)$ as given by Lemma \ref{lemma-omega}. We suppose that the conclusion of Proposition \ref{prop-psi1} holds, \emph{i.e.}
\[
\underset{\lambda\in \Lambda}{\lim_{\lambda \to 0}} \;\lambda^{-1}(\varphi(u+\lambda z)-v)
\]
is a well defined definable function from $C_u^{\Lambda,\Omega}(U)$ to $K^{m-d}$.

Then $\psi$ is a locally $1$-Lipschitz function on a dense subset of $C_u^{\Lambda,\Omega}(U)$, and $\Gamma(\psi)$ the graph of $\psi$ is dense in $C_w^\Lambda(\Gamma(\varphi))$, where $\Gamma(\varphi)$ is the graph of $\varphi$. 
\end{proposition}

\begin{proof}
By the Jacobian property, $\psi$ is differentiable on a dense subset of $C_u^{\Lambda,\Omega}(U)$. To show that $\psi$ is locally $1$-Lipschitz, by Proposition \ref{prop-tangent-lipschitz}, it suffices to show that 
\[T_{(x_0,\psi(x_0))}(\Gamma(\psi))\subseteq\set{(x,y)\in K^d \times K^{n-d} \mid \abs{y}\leq \abs{x}}=:A
\]
for $x_0$ in a dense subset of $C_u^{\Lambda,\Omega}(U)$. By definition of $\psi$, $\Gamma(\psi)\subseteq C_w^\Lambda(\Gamma(\varphi))$. As $\dim(\Gamma(\psi))=\dim(C_w^\Lambda(\Gamma(\varphi)))$, it suffices to show that $T_y C_w^\Lambda(\Gamma(\varphi))\subseteq A$ for every $y$ in a dense subset of $C_w^\Lambda(\Gamma(\varphi))$. To do so, consider the deformation $h : \overline{\Dd(\Gamma(\varphi),w,\Lambda)}\to K$ introduced in definition \ref{def-deform}. We identify the fiber $h^{-1}(0)$ with $C_w^\Lambda(\Gamma(\varphi))$ and for $\lambda\in \Lambda$, 
\[h^{-1}(\lambda)=\set{(z,\lambda)\mid w+\lambda z \in \Gamma(\varphi)}\]
with the set
\[  \set{z\in K^m \mid w+\lambda z \in \Gamma(\varphi)}.\]
As $\varphi$ is locally $1$-Lipschitz, for all $\lambda\in \Lambda$ and $y\in h^{-1}(\lambda)$, $T_y h^{-1}(\lambda)\subseteq A$ by Proposition \ref{prop-tangent-incl}. By Theorem \ref{strat-reg}, there is a $(w_h)$-stratification, hence for all $y$ in a dense subset of $C_w^\Lambda(\Gamma(\varphi))$, there is a sequence $(y_n,\lambda_n)\to (y,0)$ such that $T_y h^{-1}(0)$ is the limit of the $T_{y_n} h^{-1}(\lambda)$, hence 
\[T_y C_w^\Lambda(\Gamma(\varphi))=T_y h^{-1}(0)\subseteq A.\]

We prove now the second point. Let $z\in C_w^\Lambda(\Gamma(\varphi))$, and $(w_n)_n\in\Gamma(\varphi)$, $(\lambda_n)_n\in \Lambda$ be two sequences such that $w_n\to w$ and $\lambda_n(w_n-w)\to z$. Denote by $\pi$ the projection of $\Gamma(\varphi)$ to $U$, and $u_n=\pi(w_n)$. Then $(u_n)$ converges to $u$ and $(\lambda_n(u_n-u))$ to $a=\pi(z)\in C_u^\Lambda(U)$. We fix an $\epsilon >0$. By Lemma \ref{lemma-omega}, there is an $a'\in C_u^{\Lambda,\Omega}(U)$ with $\abs{a'-a}\leq \epsilon$ and such that $u+\lambda a'\in U$ for every $\lambda \in \Lambda$ small enough. Because $\varphi$ is globally $C$-Lipschitz, we have
\[
\abs{\lambda_n[(\varphi(\lambda_n^{-1}a'+u)-v)-(\varphi(u_n)-v)]}\leq C\abs{a'-\lambda_n(u_n-u)}. 
\]
By taking the limit, we have
\[\lim_{n\to +\infty} \abs{\lambda_n(\lambda_n^{-1}a', \varphi(\lambda_n^{-1}a'+u)-v)-\lambda_n(w_n-w)}\leq C\abs{a'-a},
\]
\[
\abs{(a',\psi(a'))-z}\leq C\epsilon, 
\]
which shows that $\Gamma(\psi)$ is dense in $C_w^\Lambda(\Gamma(\varphi))$.
\end{proof}

We can now finish the proof of Theorem \ref{thm-princ}.
\begin{proof}[Proof of Theorem \ref{thm-princ}]
We have $X\subseteq K^m$, definable of dimension less or equal to $d$. By Proposition \ref{partition-1lipsch}, there is a definable $Y$ of dimension less than $d$ and a $\kk$-definable partition of $X\setminus Y$ such that for each part $X_\xi$, there is a $\xi$-definable open $U_\xi\subseteq K^d$ and a $\xi$-definable differentiable locally 1-Lipschitz function $\phi_\xi : U_\xi \to K^{m-d}$ of graph $\Gamma_\xi$ such that for some permutation of coordinates $\gamma_\xi$, $X_\xi=\gamma_\xi(\Gamma_\xi)$. As the $\gamma_\xi$ are isometries, we have the equalities
\[
\Theta_d(\gamma_\xi(\Gamma_\xi),x)=\Theta_d(\Gamma_\xi,\gamma_xi^{-1}(x)) \mrm{and} \Theta_d(C_x^\Lambda(\gamma_\xi(\Gamma_\xi)),0)= \Theta_d(C_{\gamma_\xi^{-1}(x)}^\Lambda(\Gamma_\xi),0)
\]
By Lemma \ref{fibration}, it suffices to find a $\Lambda\in \Dd$ such that for all $\xi$, 
\[
\iota(\Theta_d(\Gamma_\xi,\gamma_\xi^{-1}(x)))=\iota(\Theta_d(C_{\gamma_\xi^{-1}(x)}^\Lambda(\Gamma_\xi),0)).
\]
Suppose $w:=\gamma_\xi^{-1}(x)=(u,v)$. By Corollary \ref{finite_accumulation}, up to refining the partition we can assume $\lim_{y\to u}\phi_\xi(y)$ exists. If this limit is not $v$, both sides of the equality we need to show are 0, so we can assume the limit is $v$. Using Remark \ref{rem-lambda-unique}, Lemmas \ref{lemma-dens1} to \ref{lemma-omega}, and Proposition \ref{prop-psi1}, we refine the partition and assume the hypotheses of Proposition \ref{prop-psi2} hold. To simplify, we keep the notations of Proposition \ref{prop-psi2}, we have a $\Lambda\in\Dd$, $U$ an open definable of $K^d$, $\varphi : U \to K^{m-d}$, $\lim_{y\to u}\varphi(y)=v$, $w=(u,v)$, $\dim(C_u^\Lambda(U))=d$, $\Omega$ such that $C_u^{\Lambda,\Omega}(U)$ is dense in $C_u^\Lambda(U)$, and $\psi : C_u^{\Lambda,\Omega}(U) \to K^{m-d}$ the definable function 
\[
\psi(z)=\underset{\lambda\in \Lambda}{\lim_{\lambda \to 0}} \;\lambda^{-1}(\varphi(u+\lambda z)-v).
\]

As $\varphi$ is locally 1-Lipschitz, $z\in U \mapsto (z,\varphi(z))$ preserves the motivic measure. By Proposition \ref{prop-tangent-incl}, for $z$ small enough, $\abs{\varphi(z)}\leq \abs{z}$, hence
\[\Theta_d(U,u)=\Theta_d(\Gamma(\varphi),w).\]
By Proposition \ref{prop-psi2}, $\psi : C_u^{\Lambda,\Omega}(U) \to K^{m-d}$ a definable locally 1-Lipschitz function on a dense subset of $C_u^{\Lambda,\Omega}$, and $\Gamma(\psi)$ is dense in $C_w^\Lambda(\Gamma(\varphi))$. By the Jacobian property, $z\mapsto (z,\psi(z))$ preserve the motivic $d$-dimensional measure and as before by Proposition \ref{prop-tangent-incl}, for $z$ small enough $\abs{\psi(z)}\leq \abs{z}$, hence
\[
\Theta_d(C_u^\Lambda(U),0)=\Theta_d(C_u^{\Lambda,\Omega}(U),0)=\Theta_d(\Gamma(\psi),0)=\Theta_d(C_w^\Lambda(\Gamma(\varphi),0)),
\]
the first and last equality begin obtained by density. Finally, by Proposition \ref{full-dim}, we have
\[
\iota(\Theta_d(U,u))=\iota(\Theta_d(C_u^\Lambda(U),0)).
\]
Putting the three equalities together, we get
\[
\iota(\Theta_d(\Gamma(\varphi)),w)=\iota(\Theta_d(C_w^\Lambda(\Gamma(\varphi)),0)).
\]
\end{proof}

\begin{remark}
\label{dimension-tan-cone}
In the above proof, we implicitly used that if $\dim(C_u^\Lambda(U))<d$, then $\dim(C_w^\Lambda(\Gamma(\varphi))<d$. If not, then there would be some coordinate projection $\pi : K^m \to K^d$ such that for $U'=\pi(\Gamma(\varphi))$, $\dim(C_{\pi(w)}^\Lambda(U'))=d$, but this would imply that $\Theta_d(\Gamma(\varphi),w)\neq 0$. 
\end{remark}

\subsection{Stabilisation of the tangent cone}

We prove in this section that the tangent cone stabilizes, that is, the existence of a distinguished tangent cone. 
\begin{theorem}
Let $X$ a definable of $K^m$ of dimension $d$. Then there is a $\Lambda\in \mathcal{D}$ such that for all $\Lambda'\in \mathcal{D}$, $\Lambda'\subseteq \Lambda$ and $x\in K^n$, we have
\[
C_x^\Lambda(X)=C_x^{\Lambda'}(X).
\]
\end{theorem}

\begin{proof}
The proof goes by induction on the dimension $d$ of $X$. If $d=0$, then $C_x^\Lambda(X)=\set{0}$ for any $\Lambda$. We suppose then that $\dim(X)=d\geq 1$ and that the proposition is true for lower dimensions. 

By Proposition \ref{partition-1lipsch}, up to taking a $\kk$-partition and working on one of the $\kk$-parts, we can suppose $X$ is a graph of a locally $1$-Lipschitz function $\varphi$ defined on an open $U \subseteq K^d$. Note also that only finitely many $\Lambda\in \Dd$ will occur, hence at the end we can pick a $\Lambda\in \Dd$ adapted to all the parts. In what follows, we fix an $x\in K^n$ to simplify the notations, but we are in fact working uniformly in $x$. 
Let $u$ be the projection of $x$ on $K^d$. Let $\Lambda\in\Dd$ be as given by Corollary \ref{cor-lambda}. Suppose first that $C_u^\Lambda(U)$ is of dimension $d$ and fix a $\Lambda'\in \Dd$, $\Lambda'\subseteq \Lambda$. Then by Corollary \ref{cor-stab-cone}, 
\[
C_u^{\Lambda'}(U)=C_u^\Lambda(U).
\]
By Propositions \ref{prop-psi1} and \ref{prop-psi2} applied to $\varphi$ and $\Lambda$, there is an $\Omega$ and a function 
\[
\psi : C_u^{\Lambda,\Omega}(U) \to K^{m-d}
\]
with graph dense in $C_x^\Lambda(X)$. These propositions also apply for $\Lambda'$, yielding an $\Omega'$ and a function 
\[
\psi' : C_u^{\Lambda',\Omega'}(U) \to K^{m-d}
\] with graph dense in $C_x^{\Lambda'}(X)$. By taking their intersections, we can suppose $\Omega=\Omega'$. Then $C_u^{\Lambda,\Omega}(U)=C_u^{\Lambda',\Omega'}(U)$ and $\psi$ and $\psi'$ agree. Hence
\[
C_x^\Lambda(X)=\overline{\Gamma(\psi)}=\overline{\Gamma(\psi')}=C_x^{\Lambda'}(X).
\]

Suppose now that $C_x^\Lambda(X)$ is of dimension $<d$ (or equivalently, $C_u^\Lambda(U)$ of dimension $<d$ by Remark \ref{dimension-tan-cone}). We will show the existence of a definable (with additional parameters) $Y\subseteq X$ such that $\dim(Y)=\dim(C_x^\Lambda(X))$ and $C_x^\Lambda(X)=C_x^\Lambda(Y)$. If we do so, we are done by induction. 

Apply the parametrized curve selection \ref{CSL-par} to $\Dd(X,x,\Lambda)\times C_x^\Lambda(X)$, with for every $y\in C_x^\Lambda(X)$, 
\[a_y=y\in \overline{\Dd(X,x,\Lambda)}\backslash \Dd(X,x,\Lambda).\]
Hence there is a nonempty definable $D$, and a definable function
\[
\sigma : B_e\times D\times C_x^\Lambda(X)
\]
such that for any $(d_0,y)\in D\times C_x^\Lambda(X)$, $\sigma_{d_0,y}(B_e^*)\subseteq\Dd(X,x,\Lambda)$ and $\lim_{t\to 0} \sigma_{d_0,y}(t)=\sigma_{d_0,y}(0)=y$. 

Then if we fix a $d\in D$, the $d$-definable set $\sigma(B_e^*,d,C_x^\Lambda(X))\subseteq \Dd(X,x,\Lambda)$ is of dimension $d+1$, and we can map it to $X$ using
\[
(z,\lambda)\in \Dd(X,x,\Lambda)\mapsto x+\lambda\cdot z\in X.
\]
The image $Y(d_0)$ is then a $d_0$-definable of dimension $d$, such that $C_x^\Lambda(Y(d_0))=C_x^\Lambda(X)$. 
\end{proof}

\section{Application to the $p$-adic case}
\label{section_application_padic}
\subsection{}
\label{cor-padic}
In this section we show how our results extend the work of Cluckers, Comte and Loeser \cite{cluckers_local_2012} on $p$-adic local density. 

For the remainder of the section, we fix $K$ a finite extension of $\Qq_p$. Let $\mathcal{L}'_{\mathrm{Mac}}=\set{0,1,+,-,\cdot,(P_n)_n>0}$ the language of Macintyre, where $P_n$ is a predicate for $n$-th powers in $K^\times$. Let $\mathcal{L}'_{\mathrm{an}}=\mathcal{L}'_{\mathrm{Mac}}\cup \set{^{-1},(K\set{x_1,...,x_n})n\geq 1}$ the analytic language, where $^{-1}$ is the field inverse and $K\set{x_1,...,x_n}$ the ring of restricted power series (\ie power series converging on $\mathcal{O}_K^n$). The definable of $K$ in the language $\mathcal{L}'_{\mathrm{Mac}}$ (resp. $\mathcal{L}'_{\mathrm{an}}$) are exactly the definable of $K$ in the three sorted language $\mathcal{L}_{\mathrm{high}}$ (resp. in the language $\mathcal{L}_{K,\mathrm{an}}$ of Example \ref{example_tame}.\ref{example_analytic}). For the remainder of this section, let $\mathcal{L}'$ be either $\mathcal{L}'_{\mathrm{Mac}}$ or $\mathcal{L}'_{\mathrm{an}}$ and (accordingly) $\mathcal{L}$ be either $\mathcal{L}_{\mathrm{high}}$ or $\mathcal{L}_{K,\mathrm{an}}$. We work in the mixed tame theory $\mathcal{T}=\mathrm{Th}(K)$ in the language $\mathcal{L}$. 

Note that the residue rings $R_n =\mathcal{O}_K/\mathcal{M}_K^n$ are finite, because the residue field $k_K$ of $K$ is finite, hence a $\kk$-partition of some definable of $X\subseteq K^n$ is in fact a finite partition of $X$. Counting points induces an isomorphism $\Qm(h[0,0,0])\simeq \Zz$. More generally, let $\mu_{n,K}$ be the Haar measure on $K^n$, normalized by $\mu_{n,K}(\mathcal{O}_K^n)=1$. This measure extends to a measure on definable sets of $K^m$ of dimension at most $n$, still denoted $\mu_{n,K}$. By \cite[Proposition 9.2]{cluckers_motivic_????},  if $\varphi\in \Cm_+(h[m,0,0])$ is of dimension at most $n$ and integrable, then 
\[
\mu_n(\varphi)=\mu_{n,K}(\varphi_K).
\]

In particular, the definition of motivic local density specializes to the $p$-adic case. That is, if $X\subseteq K^n$ of dimension $d$ and $x\in K^n$, the motivic local density $\Theta_d(X,x)\in \widetilde{\Cm_+}(\set{x})\simeq \Qq_+$ equals the local density of $X$ at $x$ of \cite[Definition 2.3.2]{cluckers_local_2012}. 

Our definition of tangent cone also coincides with the definitions of \cite[Section 3.3]{cluckers_local_2012}. Note also that the definition of cone with multiplicities $CM_x^\Lambda(X)$ given in \cite[Section 5.4]{cluckers_local_2012} coincides with our Definiton \ref{def-cone-mult}, because a $\kk$-partition in this case is a finite partition and of the isomorphism $\Qm(h[0,0,0])\simeq \Zz$. 

However, unlike what is claimed in \cite{cluckers_local_2012}, this does not coincide with the other definition of multiplicities on the tangent cone, given by \cite[Section 3.6]{cluckers_local_2012} and denoted $SC_x^\Lambda(X)$, as shown by the following example. 

\begin{example}
Fix $K=\Qq_p$ for $p\neq 2$, consider the cusp $X=\set{x^2=y^3}$ and $\Lambda=\set{\lambda\in K^\times \mid P_2(\lambda)}$ the subgroup of squares in $K^\times$. The group $\Lambda$ is of finite index in  $K^\times$ and $[K^\times/\Lambda]=4$.

One has $C_0^\Lambda(X)=\set{(0,y)\mid P_2(y)}$. Fix $(0,y_0)\in C_0^\Lambda(X)$, $y_0\neq 0$. We need to compute $\Theta_2(\mathcal{D}(X,0,\Lambda),(0,y_0,0))$. For $n\geq 3 \ord(y_0)$, we have

\begin{eqnarray*}
\mathcal{D}(X,0,\Lambda)\cap B((0,y_0,0),n)&=&\\
&&\hspace{-3cm}\set{(x,y,\lambda)\in K^2\times \Lambda \mid \min\set{\ord(x),\ord(y-y_0)}\geq n, x^2=\lambda y^3}.
\end{eqnarray*}
We have $\ord(\lambda)=\ord(x^2/y^3)=2 \ord(x)-3 \ord(y_0)$ and $\ac(\lambda)=\ac(x^2)/\ac(y_0)^3$. As $y_0$ is a square, $\ac(y_0)$ is a square and $\ord(y_0)$ is even, hence for any $(x,y)\in (K^\times)2$ with $\ord(x)\geq n$, $\ord(y-y_0)\geq n$, $x^2/y^3$ is a square, \ie in $\Lambda$. Hence $\mathcal{D}(X,0,\Lambda)\cap B((0,y_0,0),n)$ is a graph of a 1-Lipschitz function. Hence

\begin{eqnarray*}
\mu_{2,K}(\mathcal{D}(X,0,\Lambda)\cap B((0,y_0,0),n))&=&\\
&&\hspace{-5cm}\mu_{2,K}(\set{(x,y)\in (K^\times)^2\mid \ord(x)\geq n, \ord(y-y_0)\geq n})=p^{-2n},
\end{eqnarray*}

then $\Theta_2(\mathcal{D}(X,0,\Lambda),(0,y_0,0))=1$ when $y_0$ is a square and 

\[SC_0^\Lambda(X)(0,y_0)= \begin{cases} 4& \mbox{if } P_2(y_0)\\ 0& \mbox{otherwise } \end{cases}.\]

This differs from the correct value of the multiplicity
\[CM_0^\Lambda(X)(0,y_0)= \begin{cases} 2& \mbox{if } P_2(y_0)\\ 0& \mbox{otherwise } \end{cases}\]
that one computes by partitioning $X$ into graphs of 1-Lipschitz functions using
\[
 g : (x,y)\in X\mapsto (\ac(x),\ac(y)).
\]
\end{example}

Hence in the paper \cite{cluckers_local_2012}, each occurrence of $SC_x^\Lambda(X)$ must be replaced by $CM_x^\Lambda(X)$. With this modification, the first part of Theorem 3.6.1 of \cite{cluckers_local_2012} remains true, as it is exactly our Theorem \ref{thm-princ} applied to the theory of $K$. 

Note that it is unclear how to attach multiplicities in the case of non locally constant constructible functions. We do not know how to suitably modify the definition of $\nu_x^\Lambda$ in  \cite[Section 3.6]{cluckers_local_2012} in order to make true the second part of Theorem 3.6.1 of \cite{cluckers_local_2012}. 

Note also that one can easily correct the proof of this theorem. Instead of \cite[Corollary 5.3.8]{cluckers_local_2012} and \cite[Section 5.6]{cluckers_local_2012}, one can copy our proof of Theorem \ref{thm-princ}, but quoting \cite[Proposition 1.5.3]{cluckers_local_2012} instead of Proposition \ref{partition-1lipsch}, \cite[Lemma 5.3.3]{cluckers_local_2012} instead of Lemma \ref{lemma-omega}, \cite[Proposition 5.3.7]{cluckers_local_2012} instead of \ref{prop-psi2} and \cite[Proposition 3.4.1]{cluckers_local_2012} instead of Proposition \ref{prop-tangent-incl}. 

We mention also that Proposition 2.4.2 of \cite{cluckers_local_2012} is false in general. 
For example, take $\varphi_n=\mathbf{1}_{A_n}$, where $A_n=\set{x\in K\mid 0\leq\ord(x)\leq n}$. We have $\Theta_1(\varphi_n)(0)=0$ but $\sup(\varphi_n)=\mathbf{1}_{B^*}$, with $B^*$ the unit ball without 0 and $\Theta_1(\mathbf{1}_{B^*})(0)=1$.

It is used in \cite[Section 5.1]{cluckers_local_2012} to reduce the second part of \cite[Theorem 3.6.1]{cluckers_local_2012} to the first part. As we have seen that the second part is meaningless without a way to attach multiplicities to constructible functions, this is not a problem. It is also used in the proof of \cite[Section 5.2]{cluckers_local_2012}. Instead of this Proposition 2.4.2, one can use in this proof the monotone convergence theorem from classical measure theory, similarly to what we have done in the proof of Proposition \ref{full-dim}. 

\subsection{Uniformity in $p$.}

Since the pioneering work of Denef on the degree of Igusa's local zeta functions \cite{denef_igusa_1987}, many results on model theory of the $p$-adics fields and $p$-adic integration have been improved into some versions uniform in $p$. This started with the uniform $p$-adic cell decompositions theorems of (Denef-)Pas \cite{pas_uniform_1989}, \cite{pas_cell_1990} and was one of the motivation for the development of motivic integration. Denef and Loeser \cite{denef_loser_motivic_igusa} show existence of motivic zeta functions specializing for almost every $p$ to Igusa's zeta functions and later Cluckers and Loeser \cite{CL_exponential_motivic} prove a general specialization theorem for integrals of motivic constructible functions with parameters. 

In this line of ideas, we prove in this section a uniform version of Theorem 3.6.1 of \cite{cluckers_local_2012}. Let us recall the setting of \cite[Section 9.1]{CL_exponential_motivic}, in order to use the specialization principle.  Fix a number field $k$ with ring of integers $\mathcal{O}$. Let $\mathcal{A}_\mathcal{O}$ be the collection of all $p$-adic completions of finite field extensions of $k$ and for $N>0$, let $\mathcal{C}_{\mathcal{O},N}$ be the collection of all $K$ in $\mathcal{A}_\mathcal{O}$ such that the residue field of $K$ is greater than $N$ and set $\mathcal{C}_{\mathcal{O}}=\mathcal{C}_{\mathcal{O},1}$. For $K$ in $\mathcal{C}_{\mathcal{O}}$, we set $\mathcal{O}_K$ its valuation ring, $\mathcal{M}_K$ its maximal ideal, $k_K$ its residue field and $q_K=\card{k_K}$. If we fix a uniformizing parameter $\varpi_K$ of $\mathcal{O}_K$, there is a unique multiplicative map $\ac : K^\times \to k_K^\times$ extending the projection $\mathcal{O}_K^\times \to k_K^\times$ and sending $\varphi_K$ to 1. We use the three sorted language $\mathcal{L}_\mathcal{O}:=\mathcal{L}_{DP}(R)$, where $R$ is any subring of $\mathcal{O}[[t]]$ containing $\mathcal{O}[t]$. Any $K\in \mathcal{C}_{\mathcal{O}}$ can be seen as a $\mathcal{L}_\mathcal{O}$-structure by interpreting $t$ by $\varpi_K$. We work in the tame theory $\mathcal{T}$ of all $(0,0,0)$-fields in the language $\mathcal{L}_\mathcal{O}$ and definable subassignments are considered up to equivalence in this theory. If $X\subseteq h[n,r,s]$ is a definable subassignment, there is a $N>0$ such that for any $K\in\mathcal{C}_{\mathcal{O},N}$, any $\mathcal{L}_\mathcal{O}$-formula defining $X$ defines the same set $X_K\subseteq K^n\times {k_K}^r\times \Zz^s$ which does not depends on the choice of $\varphi_K$. This follows from logical compatness and Ax-Kochen-Er\v sov principle, see for example \cite{DL_motives_padic_2001}. 

A motivic constructible function $\varphi\in \Cm(Z)$ gives rise for some $N>0$ and any $K\in\mathcal{C}_{\mathcal{O},N}$ to some constructible function $\varphi_K : Z_K\to \Qq$ defined in the following way. Suppose first that $\varphi=[ Y_{ /Z}]\in \Qm(Z)$. Then $Y_K$ is a subset of $Z_k\times k_K^\ell$ for a $\ell\geq0$ and one set

\[
\varphi_K:
\begin{cases}
 Z_K \to \Qq&\\
x\mapsto \card{\pi_K^{-1}(x)},&
\end{cases}
\]

where $\pi : Y \to Z$ is the canonical projection.
This extends by linearity to general functions in $\Qm(Z)$. 

If $\varphi\in \Pp(Z)$, then $\varphi$ can be written in terms of $\eL$ and definable morphisms $\alpha : Z \to \Zz$. We replace formally each occurrence of $\eL$ by $q_K$ and each function $\alpha$ by $\alpha_K : Z_K \to \Zz$. This yields a function $\varphi_K : Z_K \to \Zz$. 

By tensor product, this defines a function $\varphi_K: Z_K \to \Qq$ for any $\varphi \in \Cm(Z)$, well defined for $K\in \mathcal{C}_{\mathcal{O},N}$ for some $N>0$ (depending on $\varphi$).

With this definitions, we can now state the specialization principle \cite[Theorem 9.1.4]{CL_exponential_motivic}. Fix a definable morphism $f : Z \to S$. Recall that we have a relative version of the measures, denoted $\mu_n^S$ and $\mu_{n,K}^{S_K}$, that assign measure to constructible function of dimension $n$ relatively to $S$, in order to deal with integrals with parameters. For any motivic constructible function $\varphi\in \Cm(Z)$, of dimension $n$ relatively to $S$ and integrable relatively to $S$, there is a $N>0$ such that for any $K\in\mathcal{C}_{\mathcal{O},N}$, 
\[
\mu_n^S(\varphi)_K=\mu_{n,K}^{S_K}(\varphi_K).
\] 
Using this result, we can prove the following uniformity theorem. For all $K\in\mathcal{C}_{\Zz}$, we write $\MV^K$ and $\Theta_d^K$ for the mean value at infinity and the local density on $K$ defined in \cite{cluckers_local_2012}. 

\begin{theorem}
Let $X\subseteq h[n,0,0]$ a definable subassignment of dimension $d$. Then there is a $\Lambda\in \mathcal{D}$ and a integer $N>0$ such that for all  $K\in\mathcal{C}_{\Zz,N}$, 
\[
\Theta_d(X)_K=\Theta_d^K(X_K),
\]
\[
\Theta_d(CM^\Lambda(X))_K=\Theta_d^K(CM^{\Lambda_K}(X_K)),
\]
and for every $x\in K^n$,
\[
\Theta_d^K(X_K,x)=\Theta_d^K(CM_x^{\Lambda_K}(X_K),0).
\]
\end{theorem} 

\begin{proof}
Given  the definition of $\Theta_d$, the first two equalities follow from the specialization principle and the fact that $\MV(\varphi)_K=\MV^K(\varphi_K)$ which follow from the definition. The last equality follows from Theorem \ref{thm-princ} and the first two equalities. 
\end{proof}

\bibliographystyle{abbrv}

\begin{thebibliography}{10}

\bibitem{bekka_regular_1993}
K.~Bekka.
\newblock Regular stratification of subanalytic sets.
\newblock {\em Bull. Lond. Math. Soc.}, 25(1):7--16, 1993.

\bibitem{cluckers_local_2012}
R.~Cluckers, G.~Comte, and F.~Loeser.
\newblock Local metric properties and regular stratifications of {$p$}-adic
  definable sets.
\newblock {\em Comment. Math. Helv.}, 87(4):963--1009, 2012.

\bibitem{cluckers_non-archimedean_????}
R.~Cluckers, G.~Comte, and F.~Loeser.
\newblock Non-archimedean {Y}omdin--{G}romov parametrizations and points of
  bounded height.
\newblock {\em Forum Math. Pi}, 3:e5, 60, 2015.

\bibitem{cluckers_approximations_2012}
R.~Cluckers and I.~Halupczok.
\newblock Approximations and {L}ipschitz continuity in {$p$}-adic
  semi-algebraic and subanalytic geometry.
\newblock {\em Selecta Math. (N.S.)}, 18(4):825--837, 2012.

\bibitem{cluckers_fields_2011}
R.~Cluckers and L.~Lipshitz.
\newblock Fields with analytic structure.
\newblock {\em J. Eur. Math. Soc. (JEMS)}, 13(4):1147--1223, 2011.

\bibitem{CLR_analytic_2006}
R.~Cluckers, L.~Lipshitz, and Z.~Robinson.
\newblock {Analytic cell decomposition and analytic motivic integration.}
\newblock {\em {Ann. Sci. \'Ec. Norm. Sup\'er. (4)}}, 39(4):535--568, 2006.

\bibitem{cluckers_b-minimality_2007}
R.~Cluckers and F.~Loeser.
\newblock b-minimality.
\newblock {\em J. Math. Log.}, 7(2):195--227, 2007.

\bibitem{cluckers_constructible_2008}
R.~Cluckers and F.~Loeser.
\newblock Constructible motivic functions and motivic integration.
\newblock {\em Invent. Math.}, 173(1):23--121, 2008.

\bibitem{CL_exponential_motivic}
R.~Cluckers and F.~Loeser.
\newblock Constructible exponential functions, motivic {F}ourier transform and
  transfer principle.
\newblock {\em Ann. of Math. (2)}, 171(2):1011--1065, 2010.

\bibitem{cluckers_motivic_????}
R.~Cluckers and F.~Loeser.
\newblock Motivic integration in all residue field characteristics for
  {H}enselian discretely valued fields of characteristic zero.
\newblock {\em J. Reine Angew. Math.}, 701:1--31, 2015.

\bibitem{comte_nature_2000}
G.~Comte, J.-M. Lion, and J.-P. Rolin.
\newblock Nature log-analytique du volume des sous-analytiques.
\newblock {\em Illinois J. Math.}, 44(4):884--888, 2000.

\bibitem{comte_equisingularite_2008}
G.~Comte and M.~Merle.
\newblock {\'E}quisingularit{\'e} r{\'e}elle. {II}. invariants locaux et
  conditions de r{\'e}gularit{\'e}.
\newblock {\em Ann. Sci. {\'E}c. Norm. Sup{\'e}r. (4)}, 41(2):221--269, 2008.

\bibitem{denef_igusa_1987}
J.~Denef.
\newblock On the degree of {I}gusa's local zeta function.
\newblock {\em Amer. J. Math.}, 109(6):991--1008, 1987.

\bibitem{denef_loser_motivic_igusa}
J.~Denef and F.~Loeser.
\newblock Motivic {I}gusa zeta functions.
\newblock {\em J. Algebraic Geom.}, 7(3):505--537, 1998.

\bibitem{DL_motives_padic_2001}
J.~Denef and F.~Loeser.
\newblock Definable sets, motives and {$p$}-adic integrals.
\newblock {\em J. Amer. Math. Soc.}, 14(2):429--469 (electronic), 2001.

\bibitem{denef_p-adic_1988}
J.~Denef and L.~van~den Dries.
\newblock {$p$}-adic and real subanalytic sets.
\newblock {\em Ann. of Math. (2)}, 128(1):79--138, 1988.

\bibitem{draper_intersection_1969}
R.~N. Draper.
\newblock Intersection theory in analytic geometry.
\newblock {\em Math. Ann.}, 180:175--204, 1969.

\bibitem{halupczok_non-archimedean_????}
I.~Halupczok.
\newblock Non-{A}rchimedean {W}hitney stratifications.
\newblock {\em Proc. Lond. Math. Soc. (3)}, 109(5):1304--1362, 2014.

\bibitem{henry_sur_1984}
J.~P. Henry, M.~Merle, and C.~Sabbah.
\newblock Sur la condition de thom stricte pour un morphisme analytique
  complexe.
\newblock {\em Ann. Sci. {\'E}cole Norm. Sup. (4)}, 17(2):227--268, 1984.

\bibitem{hironaka_stratification_1977}
H.~Hironaka.
\newblock Stratification and flatness.
\newblock In {\em Real and complex singularities (Proc. Ninth Nordic Summer
  School/{NAVF} Sympos. Math., Oslo, 1976)}, pages 199--265. Sijthoff and
  Noordhoff, Alphen aan den Rijn, 1977.

\bibitem{karzhemanov}
I.~{Karzhemanov}.
\newblock {On the cut-and-paste property of algebraic varieties}.
\newblock {\em ArXiv e-prints}, Nov. 2014.
\newblock {arXiv}:1411.6084.

\bibitem{kurdyka_$wsb_1994}
K.~Kurdyka and A.~Parusi{\'n}ski.
\newblock $(w_f)$-stratification of subanalytic functions and the \'lojasiewicz
  inequality.
\newblock {\em C. R. Acad. Sci., Paris, S{\'e}r. I}, 318(2):129--133, 1994.

\bibitem{kurdyka_densite_1989}
K.~Kurdyka and G.~Raby.
\newblock Densit{\'e} des ensembles sous-analytiques.
\newblock {\em Ann. Inst. Fourier (Grenoble)}, 39(3):753--771, 1989.

\bibitem{lelong_integration_1957}
P.~Lelong.
\newblock Int\'egration sur un ensemble analytique complexe.
\newblock {\em Bull. Soc. Math. France}, 85:239--262, 1957.

\bibitem{lion_densite_1998}
J.-M. Lion.
\newblock Densit\'e des ensembles semi-pfaffiens.
\newblock {\em Ann. Fac. Sci. Toulouse Math. (6)}, 7(1):87--92, 1998.

\bibitem{loi_verdier_1998}
T.~L. Loi.
\newblock Verdier and strict {Thom} stratifications in o-minimal structures.
\newblock {\em Illinois J. Math}, 42:347--356, 1998.

\bibitem{nowak_results_????}
K.~Nowak.
\newblock Algebraic geometry over henselian rank one valued fields.
\newblock {\em ArXiv e-prints}, 2015.
\newblock {arXiv}:1410.3280v13.

\bibitem{pas_uniform_1989}
J.~Pas.
\newblock Uniform {$p$}-adic cell decomposition and local zeta functions.
\newblock {\em J. Reine Angew. Math.}, 399:137--172, 1989.

\bibitem{pas_cell_1990}
J.~Pas.
\newblock Cell decomposition and local zeta functions in a tower of unramified
  extensions of a {$p$}-adic field.
\newblock {\em Proc. London Math. Soc. (3)}, 60(1):37--67, 1990.

\bibitem{scowcroft_structure_1988}
P.~Scowcroft and L.~van~den Dries.
\newblock On the structure of semialgebraic sets over $p$-adic fields.
\newblock {\em J. Symbolic Logic}, 53(4):1138--1164, 1988.

\bibitem{thie_lelong_1967}
P.~R. Thie.
\newblock The {L}elong number of a point of a complex analytic set.
\newblock {\em Math. Ann.}, 172:269--312, 1967.

\bibitem{thom_ensembles_1969}
R.~Thom.
\newblock Ensembles et morphismes stratifi{\'e}s.
\newblock {\em Bull. Amer. Math. Soc.}, 75:240--284, 1969.

\end{thebibliography}

\end{document}